\documentclass[preprint,10pt]{siamart251216}

\usepackage{amsmath,amssymb,mathabx,mathtools}

\usepackage{color}
\newcommand{\red}{\color{red}}
\newcommand{\green}{\color{green}}
\newcommand{\blue}{\color{blue}}

\usepackage{listings}
\usepackage{calc}

\usepackage{fancyvrb}

\usepackage{algorithm,algpseudocode}
\usepackage{algorithmicx}
\algrenewcommand\alglinenumber[1]{\footnotesize #1:} 
\newcommand{\algFontSize}{\footnotesize}

\usepackage{tabu}

\usepackage{multirow}

\usepackage{siunitx}

\usepackage{tikz}

\crefformat{figure}{#2#1#3}

\DeclareMathOperator*{\argmin}{arg\,min}

\newcommand{\mni}{\medskip\noindent}

\newcommand{\ECR}{{\rm ECR}}
\newcommand{\CR}{{\rm CR}}

\newcommand{\p}{\partial}

\newcommand{\Np}{{N_p}}

\newcommand{\f}[2]{\frac{#1}{#2}}

\def\ba#1\ea{\begin{align}#1\end{align}}

\def\bas#1\eas{\begin{align*}#1\end{align*}}

\def\bat#1\eat{\begin{alignat}{3}#1\end{alignat}}

\def\bats#1\eats{\begin{alignat*}{3}#1\end{alignat*}}

\newcommand{\bse}{\begin{subequations}}
\newcommand{\ese}{\end{subequations}}

\newcommand{\Dzt}{D_{0t}}
\newcommand{\Dpt}{D_{+t}}
\newcommand{\Dmt}{D_{-t}}
\newcommand{\dt}{\Delta t}
\newcommand{\dx}{\Delta x}

\newcommand{\eqdef}{\overset{{\rm def}}{=}}

\newcommand{\bv}{\mathbf{ b}}
\newcommand{\ev}{\mathbf{ e}}
\newcommand{\fv}{\mathbf{ f}}
\newcommand{\gv}{\mathbf{ g}}
\newcommand{\jv}{\mathbf{ j}}
\newcommand{\mv}{\mathbf{ m}}
\newcommand{\pv}{\mathbf{ p}}
\newcommand{\qv}{\mathbf{ q}}
\newcommand{\rv}{\mathbf{ r}}
\newcommand{\vv}{\mathbf{ v}}
\newcommand{\xv}{\mathbf{ x}}
\newcommand{\yv}{\mathbf{ y}}
\newcommand{\zv}{\mathbf{ z}}
\newcommand{\Vv}{\mathbf{ V}}

\newcommand{\half}{{1\over2}}

\newcommand{\Real}{{\mathbb R}}
\newcommand{\Integers}{{\mathbb Z}}
\newcommand{\Fmat}{{\mathbb F}}
\newcommand{\Gmat}{{\mathbb G}}
\newcommand{\Hmat}{{\mathbb H}}
\newcommand{\Mmat}{{\mathbb M}}

\newcommand{\zerov}{\mathbf{0}}

\newcommand{\Bc}{{\mathcal B}}
\newcommand{\Lc}{{\mathcal L}}
\newcommand{\Nc}{{\mathcal N}}
\newcommand{\Pc}{{\mathcal P}}
\newcommand{\Wc}{{\mathcal W}}

\newcommand{\Chiv}{\boldsymbol{\Chiv}}

\newcommand{\ssf}{\scriptscriptstyle}

\newcommand{\Nfreq}{\Nc_f}
\newcommand{\Nf}{{\Nc_f}}
\newcommand{\Tf}{{T_f}}
\newcommand{\Tfm}{{T_{f,m}}}
\newcommand{\Tfi}{{T_{f,i}}}

\newcommand{\sinc}{{\rm sinc}}
\newcommand{\cosc}{{\rm cosc}}

\newcommand{\vHat}{\hat{v}}
\newcommand{\vvHat}{\hat{\vv}}
\newcommand{\evHat}{\hat{\ev}}
\newcommand{\fHat}{\hat{f}}

\newcommand\floor[1]{\left\lfloor#1\right\rfloor}

\newcommand{\wFunc}{{{\mathrm{w}}}}
\newcommand{\wHat}{\hat{\wFunc}}
\newcommand{\WFunc}{{{\mathrm{W}}}}

\newcommand{\uHat}{\hat{u}}
\newcommand{\uvHat}{\hat{\mathbf{u}}}

\newcommand{\ACR}{{\rm ACR}}

\newcommand{\lambdaTilde}{\tilde{\lambda}}

\newcommand{\Ttilde}{\tilde{T}}

\newcommand{\alphaI}{\alpha^{\ssf I}}
\newcommand{\betaI}{\beta^{\ssf I}}

\newlength{\tfwidth}
\newlength{\tfheight}
\newlength{\tfxa}
\newlength{\tfxb}
\newlength{\tfya}
\newlength{\tfyb}

\newcommand{\trimFigNoBox}[6]{%
\setlength\fboxsep{1pt}
\setlength\fboxrule{0.0pt}
\fbox{\includegraphics[width=#2, clip, trim=#3 #4 #5 #6]{#1}}%
}

\newsavebox\figBox

\newcommand{\trimw}[6]{%
\sbox\figBox{\includegraphics{#1}}
\setlength{\tfwidth}{\the\wd\figBox}
\setlength{\tfheight}{\the\ht\figBox}
\setlength{\tfxa}{\tfwidth*\real{#3}}%
\setlength{\tfxb}{\tfwidth*\real{#4}}%
\setlength{\tfya}{\tfheight*\real{#5}}%
\setlength{\tfyb}{\tfheight*\real{#6}}%
\trimFigNoBox{#1}{#2}{\tfxa}{\tfya}{\tfxb}{\tfyb}%
}

\usepackage{xargs}
\newcommandx{\figByWidth}[9][5=0, 6=0, 7=0, 8=0,9=]{
\draw (#1,#2) node[anchor=south west,xshift=-16pt,yshift=-4pt] {\trimw{#3}{#4}{#5}{#6}{#7}{#8}};}

\headers{A Multi-Frequency WaveHoltz Algorithm}{A Multi-Frequency WaveHoltz Algorithm}

\title{
A Multi-Frequency Helmholtz Solver Based on the WaveHoltz Algorithm\thanks{Submitted to the editors \today.
\funding{Research supported by the National Science Foundation under grants DMS-2345225 and DMS-2513122; U.S. Department of Energy, Office of Science, Advanced Scientific Computing Research (ASCR), under award number DE-SC0025424, by the National Science Foundation under grant NSF DMS-2208164, NSF DMS-2436319, and Virginia Tech.}}
}

\author{
  Francis Appiah\thanks{Department of Mathematical Sciences, Rensselaer Polytechnic Institute, Troy, NY 12180, USA (\email{appiaf@rpi.edu}, \email{banksj3@rpi.edu}, \email{carric2@rpi.edu}, \email{henshw@rpi.edu}, \email{schwed@rpi.edu}).}
  \and Daniel Appel\"o\thanks{Department of Mathematics, Virginia Tech, Blacksburg, VA 24061 USA (\email{appelo@vt.edu})}.
  \and Jeffrey W. Banks\footnotemark[2] \and \hspace{2in} Cassandra Carrick\footnotemark[2] \and William D. Henshaw\footnotemark[2] \and D. W. Schwendeman\footnotemark[2]
}

\begin{document}

\setlength{\abovedisplayskip}{3pt}
\setlength{\belowdisplayskip}{4pt}
\setlength{\textfloatsep}{5pt plus 1.0pt minus 2.0pt}
\setlength{\floatsep}{5pt plus 1.0pt minus 2.0pt}
\setlength{\intextsep}{5pt plus 1.0pt minus 2.0pt}

\maketitle

\begin{abstract}
  We develop and analyze a new approach for simultaneously computing multiple solutions
to the Helmholtz equation for different frequencies and different forcing functions.
The new Multi-Frequency WaveHoltz (MFWH) algorithm is an extension of the original WaveHoltz method and both are based on
time-filtering solutions to an associated wave equation.
With MFWH, the different Helmholtz solutions
are computed simultaneously by solving a single wave equation combined with multiple time filters.
The MFWH algorithm defines a fixed-point iteration which can be accelerated
with Krylov methods such as GMRES.  We show that the convergence of the fixed-point iteration depends on the choice of filter parameters,
and that the choice of parameters corresponding to the original WaveHoltz algorithm can fail to converge for certain combinations of frequencies. To address this, we show how filter parameters can be found based on the solution of an optimization problem, which then improves the convergence behavior of the MFWH algorithm.
The wave equation can be solved efficiently with
implicit time-stepping using as few as five time-steps per period independent of the number of spatial grid points $N$.
If the implicit time-stepping equations are solved with an $O(N)$ solver, such as multigrid, then the MFWH algorithm has an $O(N)$ solution cost as $N$ increases (for a fixed set of frequencies).
High-order accurate approximations in space are used together with second-order accurate approximations in time.
We show how to remove time discretization errors so that the MFWH solutions converge to the corresponding solutions of the discretized Helmholtz problems.
Numerical results are given using second-order and fourth-order accurate finite-difference discretizations to confirm the convergence theory.

\end{abstract}



\section{Introduction} \label{sec:introduction}

Developing efficient solvers for the Helmholtz equation is an important and challenging topic with wide applicability in the applied sciences.
Two main challenges in solving Helmholtz equations numerically are the resolution requirements that arise
from pollution errors~\cite{BaylissGoldsteinTurkel1985,IlenburgBabuska1995,overHoltzPartTwo}, and the highly indefinite character of the discretized system of equations
which can cause significant difficulties for iterative methods~\cite{ernst2012difficult}.
There are many direct and iterative algorithms that have been developed to solve Helmholtz problems and we refer the reader to~\cite{erlangga2008advances,ModernSolversForHelmholtz2017} and the references therein for further information.

It is often the case in practice that one desires frequency domain solutions for many different frequencies or different forcings.
In this article we develop and analyze a new approach for simultaneously computing multiple solutions to the Helmholtz
equation. The scheme is an extension of the WaveHoltz algorithm which solves a single-frequency Helmholtz equation 
by time-filtering solutions to an associated wave equation~\cite{WaveHoltz1,peng2021emwaveholtz,ElWaveHoltz,WaveHoltzSemi,overHoltzPartOne}. 
For the case of multiple frequencies and forcing functions, the different Helmholtz solutions are computed simultaneously by iteratively solving a single wave equation in combination with multiple time filters.
The new scheme, called the Multi-Frequency WaveHoltz (MFWH) algorithm, defines a fixed-point iteration which can be
accelerated with Krylov methods such as GMRES.  An important element of the MFWH algorithm is the time-filtering of successive solutions of the wave equation.  The filter involves a set of parameters and the ones corresponding to the original (traditional) filter used for single-frequency problems can lead to convergence issues for the multi-frequency case.  For the fixed-point iteration, the traditional filter parameters can result in an iteration that diverges for certain choices of frequencies, while the Krylov matrix for GMRES can become singular leading to convergence to non-solutions of the Helmholtz problems.  We show how filter parameters can be chosen based on the solution of an optimization problem that restores  convergence of the fixed-point iteration and leads to nonsingular Krylov matrices.

The wave equation can be solved efficiently either with explicit time-stepping, or with
implicit time-stepping using as few as five time-steps per period independent of the number of spatial grid points $N$.
The approach we develop here has the advantage that it avoids the need to invert a large indefinite matrix for a shifted Laplacian; the implicit time-stepping matrix is definite and amenable to fast solution methods. 
When combined with an $O(N)$ solver for the implicit equations, such as multigrid, the MFWH algorithm has an $O(N)$ solution cost as $N$ increases (with the frequencies held fixed).  This $O(N)$ scaling was also shown for the single-frequency case in~\cite{overHoltzPartOne}.

Resolution requirements to overcome pollution errors are addressed through the use of high-order accurate methods.  For the purposes of this paper, we show results using second-order and fourth-order accurate discretizations, although there is no limitation to the use of even higher-order accurate approximations.  The explicit and implicit time-stepping methods for the wave equation are formally second-order accurate, and we show how time corrections can be used so that converged solutions of the MFWH algorithm agree  with the corresponding discretizations of the Helmholtz problems.  To verify the algorithm and show its convergence behavior, we focus our attention on simple geometries using finite difference discretizations on a uniform structured grid.  We also consider energy-conserving problems on closed domains using Dirichlet or Neumann boundary conditions.
The extension to multiple Helmholtz problems on complex geometries, and 
for open domains where radiation or absorbing boundary conditions are used for the wave equation, is beyond the scope this paper and is left to future work.

For reference, Table~\ref{tab:notation} provides a summary of some of the symbols and notation that are introduced in subsequent sections. The entries in the table are listed in the approximate order of appearance in the text of the paper.

\begin{table}\footnotesize
\begin{center}
  \begin{tabular}{lcl} 
     \hline 
    $\Nf$     &:& Number of different Helmholtz solutions to compute. \\
    $\omega_m$ &:& Frequencies in the Helmholtz equations, $\omega_1<\omega_2 < \ldots<\omega_{\Nf}$. \\
    $f^{(m)}(\xv)$   &:& Forcing functions for Helmholtz problems, $m=1,2,\ldots,\Nf$. \\
    $g^{(m)}(\xv)$   &:& Boundary forcing for Helmholtz problems, $m=1,2,\ldots,\Nf$. \\
    $u^{(m)}(\xv)$ &:& Continuous Helmholtz solutions, $m=1,2,\ldots,\Nf$. \\
    $\Lc$      &:& Elliptic operator in Helmholtz equation $\Lc u^{(m)} + \omega_m^2 u^{(m)} = f^{(m)}$. \\
    $\Bc$      &:& Boundary condition operator in Helmholtz equation $\Bc u^{(m)} = g^{(m)}$. \\
    $T_m$   &:& Period corresponding to $\omega_m$, $T_m=2\pi/\omega_m$. \\
    $N_{p,m}$    &:& Number of time periods for the filter of solution $m$. \\
    $T_{f,m}$ &:& Final time for the filter of solution $m$, $T_{f,m} = N_{p,m} T_m$. \\
    $\wFunc(\xv,t)$ &:& Continuous solution of the wave equation, $t\in[0,\Tf]$, $\Tf=T_{f,1}$. \\
    {$F(\xv,t)$} &:& {Composite forcing, $F(\xv,t) = \sum_{m=1}^{\Nf} f^{(m)}(\xv)\cos(\omega_m t)$.} \\
    {$G(\xv,t)$} &:& {Composite boundary forcing, $G(\xv,t) = \sum_{m=1}^{\Nf} g^{(m)}(\xv)\cos(\omega_m t)$.} \\
    $\mathcal{K}(\tau;\pv)$   &:& Filter kernel depending on a time $\tau$ and parameters $\pv$. \\
    $\pv$   &:& Filter parameters $\pv=(a,c,s)$. \\
    $\beta_m(\lambda)$ &:& Component WaveHoltz filter function for $\omega_m$. \\
    {$\beta_{d,m}(\lambda)$} &:& {Discrete (quadrature) approximation of $\beta_m(\lambda)$.} \\
    {$\tilde\beta_{d,m}(\lambda)$} &:& {Version of $\beta_{d,m}(\lambda)$ using the least-squares parameters $\tilde\pv_m$.} \\
    $v^{(m,k)}(\xv)$  &:& Continuous approximation to $u^{(m)}$ at MFWH iteration $k=0,1,2,\ldots$ \\
    $B_\omega$, $\qv$   &:& Matrix and right-hand side vector in MFWH iteration update. \\
    $(\lambda_\nu,\phi_\nu(\xv))$ &:& Eigenvalues and eigenfunctions of $(\Lc,\Bc)$, $\nu=0,1,2,\ldots$ \\
    {$\uHat_\nu^{(m)}$} &:& {Coefficient in the eigenfunction expansion of $u^{(m)}$, $\nu=0,1,2,\ldots$} \\
    {$\fHat_\nu^{(m)}$} &:& {Coefficient in the eigenfunction expansion of $f^{(m)}$, $\nu=0,1,2,\ldots$} \\
    {$\vHat_\nu^{(m,k)}$} &:& {Coefficient in the eigenfunction expansion of $v^{(m,k)}$, $\nu=0,1,2,\ldots$} \\
    {$\wHat_{\nu}$} &:& {Coefficient in the eigenfunction expansion of $\wFunc$, $\nu=0,1,2,\ldots$}\\
    $C_{\lambda_\nu}$ &:& Rank-one filter function matrix for eigenvalue $\lambda_\nu$. \\
    $M_{\lambda_\nu}$ &:& MFWH fixed-point iteration matrix, $M_{\lambda_\nu}=B_\omega\sp{-1}C_{\lambda_\nu}$. \\
    $\mu(\lambda_\nu)$ &:& Only non-zero eigenvalue of $M_{\lambda_\nu}$. \\
    $\tilde\pv$   &:& Filter parameters $\tilde\pv=(\tilde a,\tilde c,\tilde s)$ from a least-squares optimization. \\
    $L_{ph}$   &:& $p$-th order accurate approximation to $\Lc$. \\
    $B_{ph}$   &:& $p$-th order accurate boundary conditions corresponding to $\Bc$. \\
    $V_\jv^{(m,k)}$  &:& Discrete approximation of the $m$-th Helmholtz solution, $V_\jv^{(m,k)}\approx v^{(m,k)}(\xv)$. \\
    $(\alphaI,\betaI)$ &:& Time-stepping scheme weightings for the wave equation. \\
    {$F_\jv^n$} &:& Composite discrete forcing, \\
            && $F_\jv^n \eqdef \sum_{m=1}^{\Nf} f^{(m)}(\xv_\jv)\cos(\tilde\omega_m t^n)(\betaI+2\alphaI\cos(\tilde\omega_m\dt))$. \\
    {$G_\jv^n$} &:& Composite discrete boundary forcing, $G_\jv^n \eqdef \sum_{m=1}^{\Nf} g^{(m)}(\xv_\jv)\cos(\tilde\omega_m t^n)$. \\
    $\WFunc_\jv^n$  &:& Discrete solution of the wave equation $\WFunc_\jv^n\approx \wFunc(\xv_\jv,t^n)$. \\
    $N_t$ &:& Total number of time-steps in the numerical integration of the wave equation. \\
    $\tilde\omega_m$ &:& Adjusted frequencies to account for time discretization errors. \\
    {$\Ttilde_{f,m}$} &:& {Adjusted final time to account for time discretization errors.} \\
    $\mu_d(\lambda)$ &:& Discrete approximation of $\mu(\lambda)$. \\
    $\Vv_h\sp{(k)}$ &:& Vector holding all grid functions $V_\jv^{(m,k)}$, $m=1,2,\ldots,\Nf$. \\
    $A_h$ &:& Krylov matrix for the system $A_h\Vv_h=\bv_h$. \\
    $\tilde\lambda(\lambda) $ &:& Adjusted eigenvalues to account for time discretization errors. \\
    \hline
  \end{tabular}
  \end{center}
  \caption{Nomenclature. } \label{tab:notation}
\end{table}


\section{Governing equations and problem definition} \label{sec:governing}

In many applications, one wishes to find solutions to Helmholtz boundary-value problems for many different frequencies and/or
forcing functions. In this case, solutions labeled $u^{(m)}(\xv)$, $m=1,2,\ldots,\Nf$,
are sought satisfying
\bse
\label{eq:MultipleHelmholtz}
\bat
  &  \Lc  u^{(m)} + \omega_m^2 \, u^{(m)} = f^{(m)}(\xv), \qquad&&  \xv\in\Omega, \qquad \\
  &  \Bc u^{(m)} = g^{(m)}(\xv) ,  \qquad && \xv\in\partial\Omega, 
\eat
\ese
for different frequencies $\omega_m \in\Real$,  different volume forcings $f^{(m)}$, and 
different boundary data $g^{(m)}$.
Here, $\Lc= c^2 \Delta$ for the bounded domain $\Omega$, and $\Bc$ is a boundary condition operator defined over $\partial\Omega$, which we take as Dirichlet, Neumann, or Robin.\footnote{Note that the approach described here can be extended to more general 
elliptic operators $\Lc$ and more general boundary conditions.}
Note that the operators $\Lc$ and $\Bc$ are assumed to be the same for each~$m$.
We further assume that the frequencies~$\omega_m$, $m=1,2,\ldots,\Nf$,  are positive and distinct, are not at resonance, 
and are ordered from smallest to largest,
\ba 
  0< \omega_1<\omega_2< \ldots <\omega_{\Nf}.
\ea

The Multi-Frequency WaveHoltz~(MFWH) algorithm, as described in Section~\ref{sec:algorithm}, is based on solving just one forced wave equation,
\bse
\label{eq:waveEqnMultipleForcings}
  \bat
     & \p_t^2 \wFunc = \Lc \wFunc - F(\xv,t) , \qquad&&  \xv\in\Omega, \quad t\in[0,\Tf], \\
     & \Bc \wFunc = G(\xv,t)    ,        \qquad && \xv\in\partial\Omega,  \\
     & \wFunc(\xv,0) = \wFunc_0(\xv),        \qquad&&  \xv\in\Omega,     \label{eq:waveEqnMultipleForcingsIC} \\
     & \p_t \wFunc(\xv,0) = 0,            \qquad&&  \xv\in\Omega, 
  \eat
\ese
where the composite forcing function $F$ is the sum of the component forcing functions, 
and the composite boundary data function $G$ is the sum of the component boundary data functions,
\ba\label{eq:waveEqnCompositeForcings}
   & F(\xv,t) = \sum_{m=1}^{\Nf} f^{(m)}(\xv) \, \cos( \omega_m t) , \qquad
   G(\xv,t) = \sum_{m=1}^{\Nf} g^{(m)}(\xv) \, \cos( \omega_m t) .
\ea

\mni
The forced wave equation~\eqref{eq:waveEqnMultipleForcings} 
is integrated in time over an interval of length $\Tf$ 
which depends on the length of the periods $T_m \eqdef 2\pi/\omega_m$ corresponding to each frequency $\omega_m$.
In general we take the final time $\Tf$ to be 
a given integer multiple, $N_p$, of
the longest period, $T_1$, i.e.~$T_f = N_p T_1$.
Given $\Tf$, let $N_{p,m}\in\Integers^+$ denote the maximum number of periods of length $T_m$
that can fit  into the total time interval,
\bse 
\label{eq:TbarDef}
\ba
     N_{p,m} \eqdef \floor{ \f{\Tf}{T_m}},
\ea
and let $T_{f,m}$ be defined as
\ba
   T_{f,m} \eqdef N_{p,m} \, T_m , 
\ea
\ese
where $\floor{\cdot}$ is the floor function; $\floor{x}$ is the biggest integer less than or equal to $x$.


\section{The multi-frequency WaveHoltz algorithm} \label{sec:algorithm}

In this section we describe the Multi-Frequency WaveHoltz (MFWH) algorithm.
The goal of the algorithm is to solve the $\Nf$ Helmholtz problems in~\eqref{eq:MultipleHelmholtz}
for the component solutions $u^{(m)}$, $m=1,2,\ldots,\Nf$.
The algorithm is based on solving the single forced wave equation~\eqref{eq:waveEqnMultipleForcings}.
A fixed-point iteration is defined to update approximate solutions to $u^{(m)}$. This iteration
uses a set of $\Nf$ time filters that attempt to isolate the portion of the full solution $\wFunc(\xv,t)$
that has a time behavior of $\cos(\omega_m t)$.
The MFWH algorithm is described assuming continuous functions in space and time since many of the key
ingredients of the algorithm do not depend on any particular discretization. 
See~\cite{overHoltzPartOne}, for example,
for details of discrete approximations and corrections for time-discretization errors that ensure the WaveHoltz solutions
converge to the solutions of the corresponding Helmholtz BVPs even when implicit time-stepping is used with a large time-step.

Let $v^{(m,k)}(\xv)$ be an approximation to the Helmholtz solution
$u^{(m)}(\xv)$, where $k=0,1,2,\ldots$ denotes the iteration number. 
Let $\wFunc(\xv,t)$ be the solution to the forced wave equation~\eqref{eq:waveEqnMultipleForcings}
with initial condition
	\ba
	\wFunc_0(\xv) = \wFunc(\xv,0) = \sum_{m=1}^{\Nf} v^{(m,k)}(\xv) .
	\ea
Note that
the solution to the wave equation, $\wFunc(\xv,t)$, implicitly depends on 
the component iterates $v^{(m,k)}(\xv)$, $m = 1, \ldots, \Nf$, and we sometimes write $\wFunc(\xv,t; v^{(:,k)})$ to emphasize this dependence.

\newcommand{\filter}{{\mathcal{K}}}
\newcommand{\params}{{\pv_m}}
\newcommand{\paramsN}{{\pv}}
\newcommand{\paramsOrig}{{\pv_{0}}}
\newcommand{\rhs}{{q}}
\newcommand{\rhsv}{{\qv}}

\renewcommand{\algFontSize}{\small}
\begin{algorithm}
\algFontSize 
\caption{Multi-Frequency WaveHoltz Algorithm - Fixed-Point Iteration.}
\begin{algorithmic}[1]

  \Function{\blue MFWH}{$\omega_m$,$f^{(m)}$,$g^{(m)}$, $m=1,2,\ldots,\Nf$}  
    \State Compute $\Tf$, $N_m$, and $\Tfm=N_{p,m} T_m$, from $\omega_m$ and $T_m=2\pi/\omega_m$.
        \label{alg:init}
    \State $\displaystyle \beta_i(\omega_j) = \f{2}{\Tfi} \int_{0}^{\Tfi} \!\!\filter(\omega_it;\paramsN_i) \cos(\omega_j t) \, dt$ , $i,j=1,2,\ldots,\Nf$; \Comment Eval $B_\omega$~\eqref{eq:filterA} 
      \label{alg:Aentries}
    \State $k=0$; \Comment MFWH~iteration counter.
    \State $v^{(m,k)}=0,$  ~~$m=1,2,\ldots,\Nf$;  \Comment Assign initial guesses for Helmholtz iterates 
       \label{alg:v0}
    \While{ not converged} \Comment Start MFWH~iterations.

      \State $\wFunc(\xv,0) = \sum_{m=1}^\Nf v^{(m,k)}(\xv)$; \Comment Initial conditions for wave equation solve.
      \State $\wFunc(\xv,t)$= \Call{solveWaveEquation}{$\wFunc(\xv,0)$, $F$, $G$, $\Tf$}; \hspace{-5pt} \Comment Eval $\wFunc$, $t\in[0,\Tf]$. 
      \label{alg:solveWave}
      \State $\displaystyle \rhs_m(\xv) = \f{2}{\Tfm} \int_{0}^{\Tfm} \!\!\filter(\omega_mt;\params) \wFunc(\xv,t) \, dt$;
           \Comment Eval $\rhsv = [ \rhs_m ]_{m=1}^{\Nf}$~\eqref{eq:pDef}
      \State $\vv^{(k+1)}(\xv) = B_\omega^{-1} \, \rhsv(\xv)$;
                    \Comment Solve for new iterates $\vv^{(k+1)} = [ v^{(m,k+1)} ]_{m=1}^{\Nf}$
      \State $k = k+1$;
    \EndWhile    \Comment End MFWH~iterations.
    \State $u^{(m)}(\xv) = v^{(m,k)}(\xv)$ , \quad $m=1,2,\ldots,\Nf$; \Comment Approximate Helmholtz solutions.
 \EndFunction
\end{algorithmic}
\label{alg:SSWH}
\end{algorithm}

The MFWH fixed-point iteration (FPI) is given in Algorithm~\ref{alg:SSWH}.
The scheme is initialized on line~\ref{alg:init} by computing the final time, $\Tf$, 
which is determined from the longest period, as described in Section~\ref{sec:governing}.
The entries in the matrix $B_\omega$, defined below in~\eqref{eq:waveHoltzMatrixIteration},
 are evaluated on line~\ref{alg:Aentries}.
Initial guesses for $v^{(m,0)}$ (for example, $v^{(m,0)}=0$), are assigned on line~\ref{alg:v0}.
 At each iteration,
the forced wave equation is solved on line~\ref{alg:solveWave} to give $\wFunc(\xv,t; v^{(:,k)})$ using the forcing functions $F$ and $G$ defined in~\eqref{eq:waveEqnCompositeForcings}.
The filtering stage of the MFWH~algorithm uses $\Nf$ different time integrals of this
solution to define the next iterates for $v^{(m,k+1)}$, $m=1,2,\ldots,\Nf$,
\ba
\label{MFWHIteration}
  \beta_m(\omega_m) v^{(m,k+1)}(\xv) & \eqdef 
      \f{2}{\Tfm} \, \int_0^{\Tfm} \filter(\omega_mt;\params)\,
       \Big[ \wFunc (\xv,t;\, v^{(:,k)} ) \\
        &\qquad\qquad\qquad\qquad - \sum_{j=1,j\ne m}^{\Nf} v^{(j,k+1)}(\xv) \cos(\omega_j t) \Big]  \, dt,   \nonumber
\ea
where $\Tfm$ is defined in~\eqref{eq:TbarDef}.  Here, the filter functions $\beta_m(\lambda)$, $m=1,2,\ldots,\Nf$ are defined in terms of the single-frequency filter function $\beta(\lambda; \omega, T,\paramsN)$ by
\bse
 \label{eq:filterFunctionj}
\ba
  & \beta_m(\lambda) \eqdef \beta(\lambda; \omega_m, \Tfm,\params) , \\
  & \beta(\lambda; \omega, T, \paramsN) \eqdef \f{2}{T} \int_0^{T} 
      \filter(\omega t;\paramsN)\, \cos(\lambda t) \, dt ,
\ea
\ese
and the kernel $\filter$ of the filter depends on three parameters $\paramsN=(a,c,s)$ and is defined by
\ba
\label{newFilterKernel}
\filter(\tau;\,\paramsN)\eqdef a+c\,\cos(\tau)+s\,\sin(\tau).
\ea
The parameters in~\eqref{newFilterKernel} for the traditional filter are given as $\paramsOrig=(-1/4,1,0)$, but here are kept free to allow us to potentially modify the behavior of the multi-frequency filter (to be discussed).

The filter~\eqref{MFWHIteration} reduces to the usual WaveHoltz filtering step when $\Nf=1$ and $\params=\paramsOrig$.
Note that the term in the large square brackets in~\eqref{MFWHIteration} provides an approximation to $v^{(m,k)}$, in that
\ba
  \label{eq:vjApprox}
   v^{(m,k)}(\xv)\cos(\omega_mt) \approx \wFunc(\xv,t;\, v^{(:,k)}) 
         - \sum_{j=1,j\ne m}^{\Nf} v^{(j,k+1)}(\xv) \cos(\omega_j t),
\ea
and so the filter step~\eqref{MFWHIteration} approximates the usual single-frequency filter step for each individual $v^{(m,k)}$.
Also note the use of $v^{(j,k+1)}$ in the sum on the right-hand side of~\eqref{MFWHIteration}.
We call this the \textit{implicit filter} since it 
depends on the next iterate.\footnote{An explicit filter, using $v^{(j,k)}$, instead of $v^{(j,k+1)}$,
is also possible although the implicit filter appears to generally lead to faster convergence rates from our experience.}
The implicit nature of the filter is easily handled, however, as will now be shown. 
The terms in the final sum in~\eqref{MFWHIteration} can be integrated in time, leading to 
\ba
   \label{eq:filterII}
   \beta_m(\omega_m)v^{(m,k+1)}(\xv) & =         
      \f{2}{\Tfm} \, \int_0^{\Tfm} 
      \filter(\omega_mt;\params)\,  \wFunc(\xv,t;\, v^{(:,k)}) \, dt \\
       &\qquad\quad - \sum_{j=1,j\ne m}^{\Nf} \beta_m(\omega_j) \, v^{(j,k+1)}(\xv)   ,  \nonumber
\ea
and then bringing the sum on the right-hand side of~\eqref{eq:filterII} to the left-hand side gives a system of
equations for the new iterates $v^{(m,k+1)}$,
\bse
\label{eq:vkUpdate}
\ba
   & \sum_{j=1}^\Nf \beta_m(\omega_j) \, v^{(j,k+1)}(\xv)
       =  \rhs_m(\xv;v^{(:,k)}) ,     \quad m=1,2,\ldots,\Nf ,  \\
   & \rhs_m(\xv;v^{(:,k)}) \eqdef     \f{2}{\Tfm} \, \int_0^{\Tfm} \filter(\omega_mt;\params)\,  \wFunc(\xv,t;\, v^{(:,k)}) \, dt .  \label{eq:pDef}
\ea
\ese

Let $\vv^{(k)}=[v^{(1,k)}, v^{(2,k)}, \ldots, v^{(\Nf,k)} ]^T$ 
denote the vector of iterates at iteration~$k$.
Then the MFWH~filter 
stage updates~\eqref{eq:vkUpdate} can be written compactly as 
\bse
\ba
   &  B_\omega \vv^{(k+1)}  = \rhsv ,  \label{eq:waveHoltzMatrixIteration}
\ea
where the matrix $B_\omega$ is given by
\ba
 & B_\omega \eqdef 
  \begin{bmatrix}
     \beta_1(\omega_1)    & \beta_1(\omega_2)   & \cdots & \beta_1(\omega_\Nf) \\
     \beta_2(\omega_1)    & \beta_2(\omega_2)   & \cdots & \beta_2(\omega_\Nf) \\
         \vdots           &   \vdots            &        & \vdots              \\
      \beta_\Nf(\omega_1) & \beta_\Nf(\omega_2) & \cdots & \beta_\Nf(\omega_\Nf) 
  \end{bmatrix} ,  \label{eq:filterA} 
\ea
\ese
and the right-hand-side vector function $\rhsv$ has components $\rhs_m$ from~\cref{eq:pDef}.
Note that to simplify the presentation, the wave equation solver in Algorithm~\ref{alg:SSWH} returns the solution
at all grid points and all time-steps. In practice this should normally be avoided to save memory usage
and instead the intermediate variables $\rhs_m(\xv;v^{(:,k)}) $ should be accumulated as
the time-stepping progresses.

As in the single-frequency WaveHoltz algorithm~\cite{WaveHoltz1,overHoltzPartOne} 
the MFWH update~\cref{eq:waveHoltzMatrixIteration} can be written using an affine operator $\Wc$ and linear operator $S$,  
\ba
    \vv^{(k+1)} \eqdef \Wc \, \vv^{(k)} = S \, \vv^{(k)} + \bv = B_\omega^{-1} \rhsv  , \label{eq:waveHoltzMapping}
\ea
which maps $\vv^{(k)}$ to $\vv^{(k+1)}$ by solving the wave equation with initial condition $\vv^{(k)}$ and then filtering in time.

\section{Analysis of the multi-frequency WaveHoltz algorithm} \label{sec:analysis}

An analysis of the MFWH~algorithm is performed
using eigenfunction expansions to determine
the asymptotic convergence rate of the MFWH fixed-point iteration (FPI).
The theory is confirmed by numerical computations in Section~\ref{sec:numerical_confirmation_of_the_convergence_analysis}.

\subsection{Fixed-point iteration convergence rate from an eigen\-function analysis}  \label{sec:eigenFunctionAnalysis}

\newcommand{\Bmat}{{C}}

Let us assume that the eigenvalue boundary-value problem (BVP) given by 
\bse
\label{eq:eigBVP}
\bat
  &  \Lc \phi = - \lambda^2 \, \phi, \qquad&&  \xv\in\Omega, \qquad \\
  &  \Bc \phi = 0,                   \qquad && \xv\in\partial\Omega, 
\eat
\ese
has a complete set of linearly independent eigenfunctions, $\phi_\nu(\xv)$, with corresponding eigenvalues $\lambda_\nu \ge 0$, for $\nu=0,1,2,\ldots$.
This imposes some restrictions on the coefficients in any Robin boundary condition~\cite{StraussBook2007}.
Consider first the individual Helmholtz problems in~\eqref{eq:MultipleHelmholtz} with the boundary forcings $g^{(m)}(\xv)$ in~\eqref{eq:MultipleHelmholtz} set to zero to simplify the analysis.
Let $u^{(m)}(\xv)$ and $f^{(m)}(\xv)$ be expanded in eigenfunction expansions
\ba
  & u^{(m)}(\xv) = \sum_{\nu=0}^\infty \uHat^{(m)}_{\nu} \phi_\nu(\xv), \qquad
    f^{(m)}(\xv) = \sum_{\nu=0}^\infty \fHat^{(m)}_{\nu} \phi_\nu(\xv), \label{eq:eigenExpansion}
\ea
where $\uHat^{(m)}_{\nu}$ and $\fHat^{(m)}_{\nu}$ denote generalized Fourier coefficients.
Substituting~\eqref{eq:eigenExpansion} into~\eqref{eq:MultipleHelmholtz} gives expressions
for the coefficients in the expansion of the Helmholtz solution
\ba
  \uHat^{(m)}_{\nu} = \frac{ \fHat^{(m)}_{\nu} }{\omega_m^2 - \lambda_\nu^2},\qquad \nu=0,1,2,\ldots,
\ea
assuming $\omega_m \ne \lambda_\nu$, i.e.~$\omega_m$ is not at resonance.
 
To study the MFWH~iteration we first consider an eigenfunction expansion of the solution of the wave equation of the form
\bse
\ba
   \label{eq:wHatEig}
   \wFunc(\xv,t)   = \sum_{\nu=0}^\infty \wHat_{\nu}(t) \, \phi_\nu(\xv) .
\ea
Substituting this expansion, and the one in~\eqref{eq:eigenExpansion} for $f^{(m)}(\xv)$, into the wave equation IBVP~\eqref{eq:waveEqnMultipleForcings}
leads to an ODE for each coefficient $\wHat_\nu(t)$ of the form
\ba
  &  \p_t^2 \wHat_\nu = -\lambda_\nu^2  \wHat_\nu  - \sum_{j=1}^{\Nf} \fHat^{(j)}_{\nu} \, \cos(\omega_j t) , \\
  \intertext{with initial conditions}
  &  \wHat_\nu(0) = \sum_{j=1}^{\Nf} \vHat^{(j,k)}_{\nu} , \qquad
     \p_t \wHat_\nu(0) =0 ,
\ea
where $\vHat^{(j,k)}_{\nu}$ are generalized Fourier coefficients in the expansion
\ba
   \label{eq:vHatEig}
   v^{(m,k)}(\xv) = \sum_{\nu=0}^\infty \vHat^{(m,k)}_{\nu} \, \phi_\nu(\xv) .
\ea
The solution to each ODE is 
\ba
  \wHat_\nu(t) =   \sum_{j=1}^{\Nf} \left\{ \big(\vHat^{(j,k)}_{\nu} - \uHat^{(j)}_{\nu} \big)\, \cos(\lambda_\nu t) 
                       +                \uHat^{(j)}_{\nu} \cos(\omega_j t) \right\} .
      \label{eq:wHatFormula}
\ea
\ese
We now substitute the eigenfunction expansions in~\eqref{eq:wHatEig} and~\eqref{eq:vHatEig}, along with
the expression for $\wHat_\nu(t)$ in~\eqref{eq:wHatFormula},
into the MFWH~filter formula~\eqref{eq:vkUpdate} to obtain
\ba
\label{eq:solutionIterates}
  \sum_{j=1}^{\Nf} \beta_m(\omega_j) \, \vHat_\nu^{(j,k+1)}
      =  \sum_{j=1}^{\Nf} \left\{ \big(\vHat^{(j,k)}_{\nu} - \uHat^{(j)}_{\nu} \big)\, \beta_m(\lambda_\nu) 
                       +                \uHat^{(j)}_{\nu} \beta_m(\omega_j) \right\} ,
\ea
for $m=1,2,\ldots,\Nf$.

Collecting the results in~\eqref{eq:solutionIterates} as a system gives
\ba
     B_\omega \vvHat_\nu^{(k+1)} = \Bmat_{\lambda_\nu} (\vvHat_\nu^{(k)}  - \uvHat_\nu) + B_\omega \uvHat_\nu, 
\ea
where the matrix $B_\omega$ was given previously in~\eqref{eq:filterA} 
and $\Bmat_\lambda$ is a rank-one matrix defined by
\ba
 \label{eq:Bmatrix}
 & \Bmat_\lambda \eqdef 
  \begin{bmatrix}
     \beta_1(\lambda) \\
     \beta_2(\lambda) \\
         \vdots       \\
     \beta_{\Nfreq}(\lambda)      
  \end{bmatrix}
  \begin{bmatrix}
  1 & 1 & \cdots & 1
  \end{bmatrix}
  =
  \begin{bmatrix}
     \beta_1(\lambda) & \beta_1(\lambda) & \cdots & \beta_1(\lambda) \\
     \beta_2(\lambda) & \beta_2(\lambda) & \cdots & \beta_2(\lambda) \\
         \vdots       &   \vdots         &        &    \vdots        \\
     \beta_{\Nfreq}(\lambda) & \beta_{\Nfreq}(\lambda) & \cdots & \beta_{\Nfreq}(\lambda)     
  \end{bmatrix}   
    .
\ea
The equation for the vector of errors, $\evHat_\nu^{(k)}=\vvHat_\nu^{(k)}-\uvHat_\nu$, is therefore
\ba
     B_\omega \evHat_\nu^{(k+1)} = \Bmat_{\lambda_\nu} \evHat_\nu^{(k)}.
\ea
Whence the vector of errors satisfies the fixed-point iteration (FPI)
\bse
\ba
\label{ereqImp}
  \evHat^{(k+1)}_\nu = M_{\lambda_\nu}  \, \evHat^{(k)}_\nu,
\ea
where the FPI iteration matrix is 
\ba
  M_{\lambda_\nu} \eqdef B_\omega^{-1} \Bmat_{\lambda_\nu}.
\ea
\ese
Thus, the iteration converges provided the spectral radius of $M_{\lambda_\nu}$ is less than one for all $\nu$.
The iteration matrix $M_{\lambda_\nu}$ has only one non-zero eigenvalue since $\Bmat_{\lambda_\nu}$ is rank~one,
and let $\mu(\lambda_\nu)$ be this one non-zero eigenvalue. 
The asymptotic convergence rate (ACR) of the FPI, denoted by $\mu_\infty$, is thus given by 
\ba
      \mu_\infty \eqdef \sup_{\lambda_\nu} | \mu(\lambda_\nu) | . \label{eq:muInfinity}
\ea
These results are summarized in Theorem~\ref{thm:MFWaveHoltzConvergence}.

\begin{theorem} \label{thm:MFWaveHoltzConvergence}
  The MFWH~fixed-point iteration (FPI) 
  is governed by the error equation
\bse
\ba
  & \evHat^{(k+1)}_\nu = M_{\lambda_\nu} \, \evHat^{(k)}_\nu, \\
  & M_{\lambda_\nu} = B_\omega^{-1} \Bmat_{\lambda_\nu},
\ea  
\ese
for the component $\nu$ in the eigenfunction expansion.
The matrices $B_\omega$ and $\Bmat_\lambda$ are given in~\eqref{eq:filterA} and~\eqref{eq:Bmatrix}, respectively.
The matrix $M_{\lambda_\nu}$ has one non-zero eigenvalue, $\mu(\lambda_\nu)$, and the FPI converges iff 
$|\mu(\lambda_\nu)|<1$ for all eigenvalues $\lambda_\nu$.  
The asymptotic convergence rate (ACR) of the FPI is $\mu_\infty$ given in~\eqref{eq:muInfinity}.
\end{theorem}

\newcommand{\ai}{d}
The eigenvalues of $M_{\lambda_\nu}$ are roots of its characteristic polynomial $P_\Nf(z)$ which takes the form
\ba
  P_\Nf(z) \eqdef \det( z I - M_{\lambda_\nu} ) = \det( z I - B_\omega^{-1} \Bmat_{\lambda_\nu} ) = z^{\Nf} - \mu(\lambda_\nu) z^{\Nf-1} ,
\ea
since there are $\Nf-1$ zero eigenvalues. Let the entries of $B_\omega^{-1}$ be denoted by $\ai_{ij}$,
\ba
    B_\omega^{-1} = \bigl[ \ai_{ij} \bigr]. 
\ea
The following theorem shows that $\mu(\lambda_\nu)$ is a linear combination of the filter functions $\beta_m(\lambda_\nu)$.
\newcommand{\weight}{{\gamma}}
\begin{theorem} \label{thm:MFWHmu}
The non-zero eigenvalue of $M_{\lambda_\nu}$ is given by 
\ba
    \mu(\lambda_\nu) = \sum_{m=1}^{\Nf} \weight_m \, \beta_m(\lambda_\nu) ,
    \label{eq:MFWHmu}
\ea
where $\weight_m$ is the sum of the entries in column $m$ of $B_\omega^{-1}$, 
\ba
    \weight_m = \sum_{i=1}^{\Nf} \ai_{im} .
\ea
\end{theorem}
\begin{proof}
The matrix $M_{\lambda_\nu} = B_\omega^{-1} \Bmat_{\lambda_\nu}$ has entries $m_{ij}$ given by
\ba
   m_{ij} = \sum_{k=1}^{\Nf} \ai_{ik} \,\beta_k(\lambda_\nu) \eqdef m_i ,
\ea
since $m_{ij}$ does not depend on $j$. Whence
\ba
  \det( z I - M_{\lambda_\nu} ) & = 
 \det\begin{bmatrix*}[r]
     z - m_1   &  -m_1  &  -m_1 & \ldots \\
        -m_2   & z-m_2  &  -m_2 & \ldots \\
        -m_3   &  -m_3  & z-m_3 & \ldots \\
        \vdots\;\; & \vdots\;\; & \vdots\;\; &  \ddots  \\
     \end{bmatrix*} .
  \label{eq:Mdet}
\ea
From the usual Leibniz formula for the determinant, it is easy
to see that, as a function of $z$, the leading term is $z^{\Nf}$ while the coefficient 
of $z^{\Nf-1}$ is just the negative of the trace of $M_{\lambda_\nu}$, i.e.~$-{\rm Tr}(M_{\lambda_\nu})=-(m_1+m_2+\ldots +m_{\Nf})$. Whence
\ba
  \det( z I - M_{\lambda_\nu} )  = z^{\Nf} - (m_1+m_2 +\ldots +m_{\Nf}) z^{\Nf-1} .
\ea
Thus the non-zero root is
\ba
  \mu(\lambda_\nu)
      &= \sum_{i=1}^{\Nf} m_i =  \sum_{i=1}^{\Nf} \sum_{k=1}^{\Nf} \ai_{ik} \, \beta_k(\lambda_\nu)  
       = \sum_{k=1}^{\Nf} \Big[ \sum_{i=1}^{\Nf}  \ai_{ik} \Big] \, \beta_k(\lambda_\nu)      
       = \sum_{k=1}^{\Nf} \weight_k \, \beta_k(\lambda_\nu) \,,
\ea
which completes the proof.
\end{proof}

We have assumed that none of the frequencies $\omega_m$ are at resonance, i.e., no $\omega_m$ is equal to an eigenvalue $\lambda_\nu$.
This resonance condition is reflected in $\mu(\lambda)$ being equal to one when $\lambda=\omega_m$ as shown in the next theorem. 
\begin{theorem} \label{thm:muOne}
The multi-frequency filter function $\mu(\lambda)$ is equal to one at each of the frequencies $\omega_m$, $m=1,2,\ldots,\Nf$,
\ba
    \mu(\omega_m) =1, \quad m=1,2,\ldots,\Nf.     \label{eq:muOne}
\ea
\end{theorem}
\begin{proof}
Let $\rv$ be an eigenvector of $M_\lambda$.
Since $M_\lambda= B_\omega^{-1} \Bmat_\lambda$ it follows that 
\ba
  \Bmat_\lambda \rv = \mu(\lambda) B_\omega \rv .
\ea
Let $\rv=[0,\ldots,0,1,0,\ldots,0]^T$ be the unit vector in $\Real^{\Nf}$ in direction $m$.
Then from the definition of $B_\omega$ in~\eqref{eq:filterA} and $\Bmat_\lambda$ in~\eqref{eq:Bmatrix} it follows that 
\ba
  B_\omega \rv = \Bmat_{\omega_m}\rv = 
      \begin{bmatrix}
             \beta_1(\omega_m) \\
             \beta_2(\omega_m) \\
             \vdots \\
             \beta_\Nf(\omega_m) 
          \end{bmatrix} .
\ea
Thus, $M_{\omega_m} \rv = \rv$ when $\lambda=\omega_m$, and therefore $\rv$ is an eigenvector of $M_{\omega_m}$ with eigenvalue $\mu(\omega_m)=1$.
This completes the proof.
\end{proof}

\subsection{Form of the MFWH~fixed-point iteration function $\mu$}  \label{sec:formOfMu}

\begin{figure}[H]
	\centering
	\includegraphics[width=0.35\textwidth]{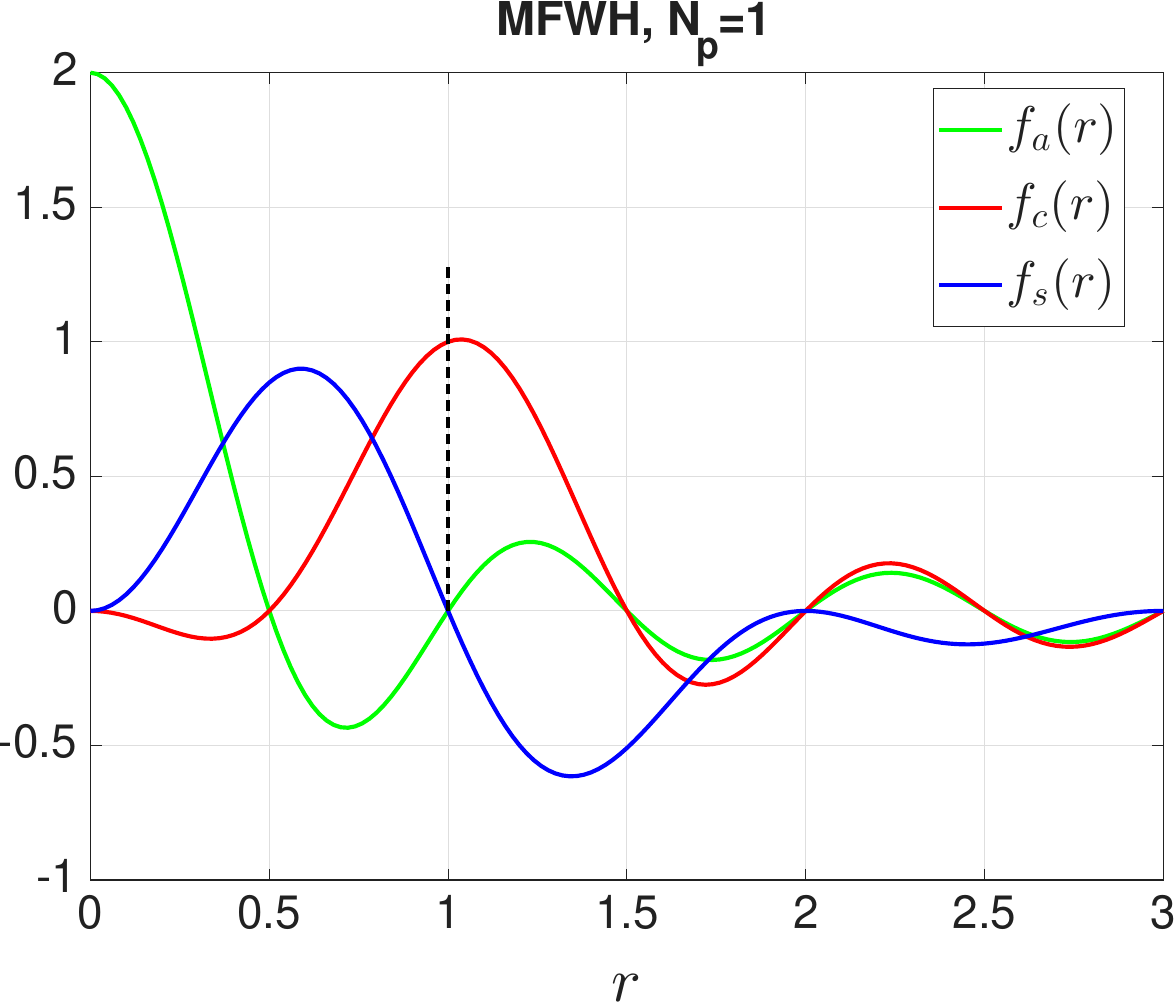}\hspace{.2in}
	\includegraphics[width=0.35\textwidth]{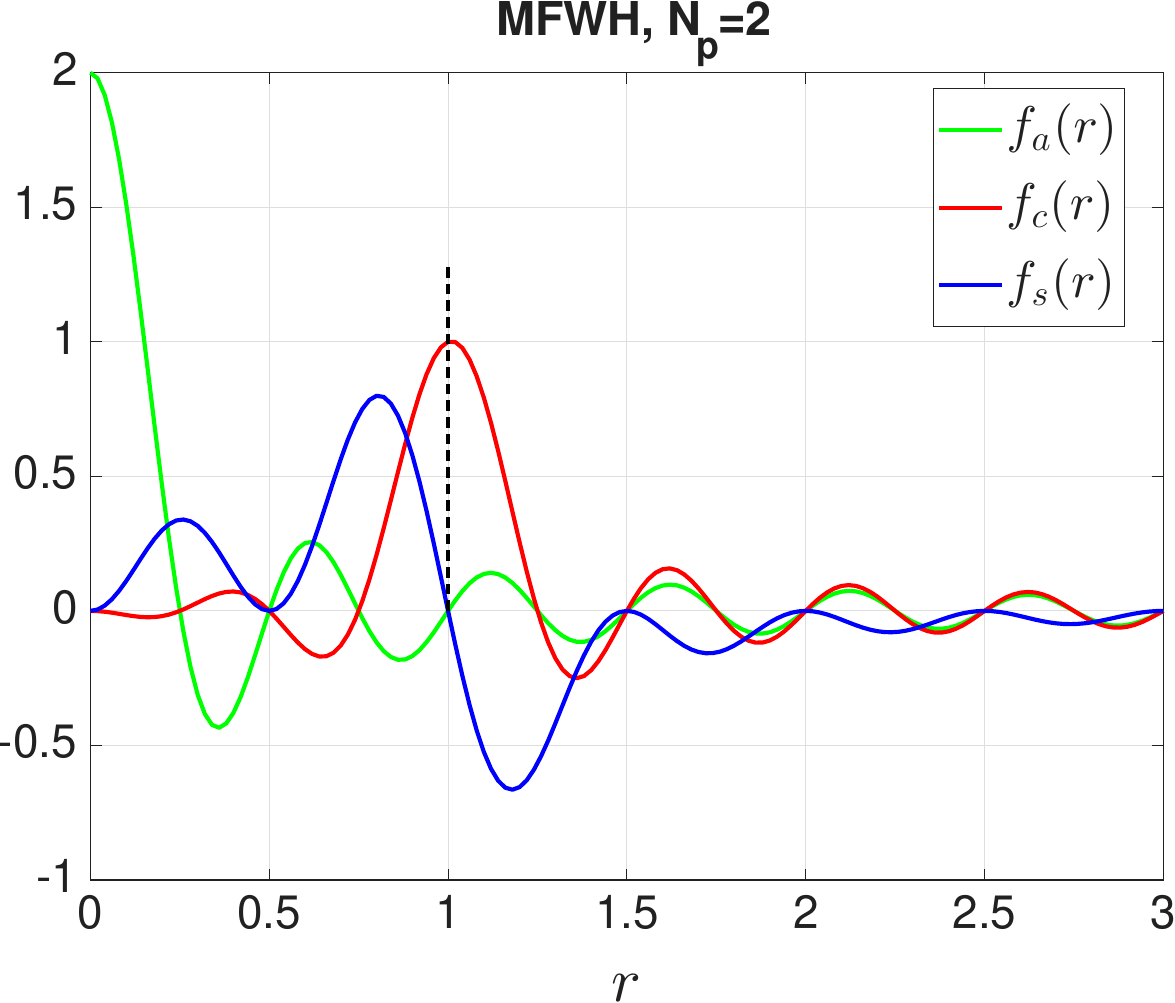}
	\caption{Plots of $f_a(r)$, $f_c(r)$, and $f_s(r)$ for $r=\lambda/\omega$ and $T = 2\pi \Np/\omega$.  Left: $\Np=1$; Right: $\Np=2$.}
	\label{fig:BetaTerms}
\end{figure}

The convergence of the MFWH~fixed-point iteration is governed by the function $\mu=\mu(\lambda)$ in~\eqref{eq:MFWHmu}. 
It has been established that $\mu(\lambda)$ is a weighted sum of filter functions $\beta_m(\lambda)$.
Each filter function $\beta_m(\lambda)$ is, in turn, the sum of three component functions,
\ba
\label{eq:betafct}
   \beta_m(\lambda) &= \beta(\lambda; \omega_m,T_{f,m},\params) \\
    &= a_mf_a(\lambda;\omega_m)+c_mf_c(\lambda;\omega_m,T_{f,m})+s_mf_s(\lambda;\omega_m,T_{f,m}), \nonumber
\ea
where
\bse
\label{eq:compFuncs}
\ba
		f_{a}(\lambda;T) &\eqdef \frac{2}{T}\int_{0}^{T}\cos(\lambda t)\,dt = 2\, \sinc(\lambda T),\\
		f_{c}(\lambda;\omega,T) &\eqdef \frac{2}{T}\int_{0}^{T}\cos(\omega t)\cos(\lambda t)\,dt = \sinc((\omega-\lambda)T) + \sinc((\omega+\lambda)T),\\
		f_{s}(\lambda;\omega,T) &\eqdef \frac{2}{T}\int_{0}^{T}\sin(\omega t)\cos(\lambda t)\,dt=\cosc((\omega-\lambda)T) + \cosc((\omega+\lambda)T),
\ea
\ese
and
\ba
\sinc(x) = \frac{\sin(x)}{x}, \qquad \cosc(x) = \frac{1-\cos(x)}{x}.
\ea
The behavior of the component functions in~\eqref{eq:compFuncs} is displayed in Figure~\ref{fig:BetaTerms}.  All three component functions oscillate with their oscillations, becoming sharper as $\Np$ increases.  Of particular interest is the behavior near $\lambda=\omega$ as this generally plays an important role in the convergence, or lack of convergence, of the fixed-point iteration since $\mu(\lambda)$ is a linear combination of $\beta_m(\lambda)$, $m=1,2,\ldots,\Nf$.  We observe that $f_c=1$ and $f_a=f_s=0$ at $\lambda=\omega$ and so the peak in $\beta_m$ at $\omega_m$ is controlled by the parameter $c_m$.  The slope of $\beta_m$ at $\omega_m$ is determined by all three parameters in~$\params$, and we note that it is often desirable to adjust the full set of parameters so that $\partial_\lambda \mu=0$ at $\lambda=\omega_m$, $m=1,2,\ldots,\Nf$, which makes these points local extrema. 

For the multi-frequency case $\mu(\lambda)$ has the form in~\eqref{eq:MFWHmu}.  With $\Nf=2$, for example, $\mu(\lambda)$ takes the form
\ba
  &  \mu(\lambda) = \weight_1 \beta_1(\lambda) +
                    \weight_2 \beta_2(\lambda), 
\ea
where the weights are given by
\bse
\ba
\weight_1&=\frac{\beta_2(\omega_2) - \beta_2(\omega_1)}{\beta_1(\omega_1)\beta_2(\omega_2)-\beta_2(\omega_1)\beta_1(\omega_2)},\\
\weight_2&=\frac{\beta_1(\omega_1) - \beta_1(\omega_2)}{\beta_1(\omega_1)\beta_2(\omega_2)-\beta_2(\omega_1)\beta_1(\omega_2)}.
\ea
\ese

{
\newcommand{\figSize}{5.3cm}
\begin{figure}[htb]
\begin{center}
\resizebox{\textwidth}{!}{
\begin{tikzpicture}[scale=1]
  \useasboundingbox (0,.35) rectangle (3*\figSize,4.5);  

   \begin{scope}[yshift=0cm]
    \figByWidth{0*\figSize}{0}{fig/OmfwhOmega5Omega9Np1BddddFixPointMuFunction}{\figSize}[.0][0.0][0.][0.]
    \figByWidth{1*\figSize}{0}{fig/OmfwhOmega5Omega9Np2BddddFixPointMuFunction}{\figSize}[.0][0.0][0.][0.]
    \figByWidth{2*\figSize}{0}{fig/OmfwhOmega5Omega9Np3BddddFixPointMuFunction}{\figSize}[.0][0.0][0.][0.] 
       
  \end{scope}  
\end{tikzpicture}
}
\end{center}
\caption{Multi-frequency WaveHoltz filter function for two frequencies,  $\omega_1=5$ and
$\omega_2=9$. The solution is integrated over $\Np$ periods of the longest period $T_1=2\pi/\omega_1$.
    } 
\label{fig:MFWHfilterFunction2}
\end{figure}
}

Figure~\ref{fig:MFWHfilterFunction2} plots $\mu(\lambda)$ for a two-frequency case with $\omega_1=5$ and
$\omega_2=9$ with filter parameters taken to be the traditional choice, $\params=\pv_0=(-1/4,1,0)$, $m=1,2$.  
The left graph shows $\mu$ when the number of periods (for the longest period $T_1$) is $N_p=1$.
The component filters $\beta_m$ are also graphed.
 Note that $\beta_2$ is quite wide since $N_{p,2}=1$ (i.e.~only one period of mode 2 fits into the time interval $T_f=T_1$).
The middle and right graphs show $\mu$ for $N_p=2$ and $N_p=3$, respectively, and we observe that the behavior of $\mu$ is more sharply peaked near $\omega_1$ and $\omega_2$.
For $N_p=1$ and $N_p=2$, there are intervals where $|\mu|>1$, and the FPI would not
converge if an eigenvalue $\lambda_\nu$ lies in these intervals. 
Although a problem for the FPI, this is generally not a problem
for the GMRES accelerated algorithm discussed in Section~\ref{sec:algorithm_gmres}.
A more subtle problem occurs if there is a spurious resonance where $\lambda_\nu\ne\omega_m$, $m=1,2,\ldots,\Nf$,  but $\mu(\lambda_\nu)=1$ (or a near spurious resonance where $\mu(\lambda_\nu)$ is very close to~$1$).  If a spurious resonance occurs, then the matrix forming the Krylov basis in GMRES becomes singular and GMRES may not converge, or it may converge to a non-solution, see the discussion in Section~\ref{sec:algorithm_gmres}.

In the next section, we describe a procedure for choosing the filter parameters so that intervals with $|\mu|>1$ are suppressed and $\vert\mu\vert$ is small in intervals away from $\lambda=\omega_m$.

\section{Filter parameters from an optimization problem} \label{sec:optimization}

We now describe a procedure for selecting the parameters $\params=(a_m,c_m,s_m)$ in the filter functions~\eqref{eq:filterFunctionj}, $m=1,2,\ldots,\Nf$.  Recalling that $\mu(\lambda)=1$ when $\lambda=\omega_m$, $m=1,2,\ldots,\Nf$, independent of the choice of the parameters, our basic objective is to select a set of parameters so that $\vert\mu(\lambda)\vert$ is small for $\lambda\ne\omega_m$, and in particular $\vert\mu(\lambda)\vert<1$ so that the FPI converges and spurious resonances are avoided. To do this, let us first define the set of parameters
\newcommand{\Pset}{{\Pc}}
\newcommand{\PsetOpt}{{{\Pc}_{{\rm opt}}}}
\ba
\Pset\eqdef\bigl\{\params\;\vert\; m=1,2,\ldots,\Nf\bigr\},
\ea
and then construct an optimization problem by setting
\ba
\label{eq:Popt}
\PsetOpt=\argmin_{\Pset} J(\Pset),
\ea
where $J$ is a suitable objective function depending on the set of parameters. Many choices for $J$ could be made, but the one chosen here is
\ba
\label{eq:Jfct}
J(\Pset)\eqdef \int_\Lambda \mu(\lambda;\,\Pset)\sp2\,d\lambda,
\ea
where $\Lambda$ is a chosen (fixed) interval  containing the frequencies $\omega_m$, $m=1,2,\ldots,\Nf$. 
The objective function in~\eqref{eq:Jfct} corresponds to an $L_2$ minimization.\footnote{An $L_1$ optimization has also been considered and the results are not significantly different.}  In addition, the optimization problem is constrained by setting
\bse
\label{eq:constraints}
\ba
c_m=1,&\qquad m=1,2,\ldots,\Nf , \label{eq:cone}\\
\bigl.\partial_\lambda\mu(\lambda;\,\Pset)\bigr\vert_{\lambda=\omega_m}=0,&\qquad m=1,2,\ldots,\Nf .\label{eq:ctwo}
\ea
\ese
The $\Nf$ constraints in~\eqref{eq:cone} are imposed to fix a scale for each filter function, $\beta_m$, used in~\eqref{MFWHIteration}, while the $\Nf$ ones in~\eqref{eq:ctwo} are necessary conditions for the points $\mu(\omega_m;\,\Pset)=1$ to be local extrema.  The constrained optimization problem described in~\eqref{eq:Popt}, \eqref{eq:Jfct} and~\eqref{eq:constraints} is nonlinear and involves $\Nf$ free parameters.
While a solution procedure for this problem could be pursued,
solving it when $\Nf$ is large could be computationally challenging.\footnote{This approach was initially pursued but found to be expensive, even for moderately small values of $\Nf$.}

Fortuitously, the form of the filter function enables the problem to be recast as an unconstrained linear optimization problem, which can then be solved readily using a least-squares approach.  To show this, consider first the iteration step in~\eqref{MFWHIteration} and the definition of $\beta_m$ in~\eqref{eq:filterFunctionj} and~\eqref{newFilterKernel}.  Multiplication of~\eqref{MFWHIteration} by the weight $\weight_m$ from the formula~\eqref{eq:MFWHmu} for $\mu(\lambda)$ gives
\newcommand{\betaT}{{\tilde\beta}}
\newcommand{\paramsT}{{\tilde\paramsN}}
\ba
\label{MFWHIterationT}
  \betaT_m(\omega_m) v^{(m,k+1)}(\xv) & \eqdef
      \f{2}{\Tfm} \, \int_0^{\Tfm} \filter(\omega_mt;\paramsT_m)\,
       \Big[ \wFunc (\xv,t;\, v^{(:,k)} ) \\
        &\qquad\qquad\qquad\qquad - \sum_{j=1,j\ne m}^{\Nf} v^{(j,k+1)}(\xv) \cos(\omega_j t) \Big]  \, dt,   \color{black} \nonumber
\ea
where
\ba
\betaT_m(\lambda) \eqdef \beta(\lambda; \omega_m, \Tfm,\paramsT_m),\qquad \paramsT_m=\weight_m\params.
\ea
Note that the form of the iteration step in~\eqref{MFWHIterationT} is unchanged from the original one, and that the formula for $\mu(\lambda)$ in~\eqref{eq:MFWHmu} now reads
\ba
\mu(\lambda)&=\sum_{m=1}\sp\Nf \betaT_m(\lambda)  \label{eq:linearOptMu}\\
      &= \sum_{m=1}\sp\Nf\Big\{ \tilde{a}_m f_a(\lambda;\omega_m)
                              +\tilde{c}_m f_c(\lambda;\omega_m,T_{f,m})+ \tilde{s}_m f_s(\lambda;\omega_m,T_{f,m}) \Big\}, \nonumber 
\ea
which is a linear function of the (weight-scaled) parameters $\paramsT_m$, $m=1,2,\ldots,\Nf$.  The corresponding optimization problem is
\newcommand{\PsetT}{{\tilde\Pc}}
\newcommand{\PsetOptT}{{{\tilde\Pc}_{{\rm opt}}}}
\ba
\label{eq:PoptT}
\PsetOptT=\argmin_{\PsetT} \int_\Lambda \left(\sum_{m=1}\sp\Nf \betaT_m(\lambda)\right)\sp2 \,d\lambda,
\ea
subject to the constraints
\bse
\label{eq:constraintsT}
\ba
\sum_{j=1}\sp\Nf \betaT_j(\omega_m)=1,&\qquad m=1,2,\ldots,\Nf , \label{eq:coneT}\\
\Bigl.\sum_{j=1}\sp\Nf \partial_\lambda\betaT_j(\lambda)\Bigr\vert_{\lambda=\omega_m}=0,&\qquad m=1,2,\ldots,\Nf .\label{eq:ctwoT}
\ea
\ese
Note that in the original optimization problem~\eqref{eq:Popt}, $\mu(\omega_m)=1$ for any values of the parameters $\pv_m$.
However, this is no longer true for~\eqref{eq:PoptT} and thus the constraint~\eqref{eq:coneT} replaces the normalization constraint~\eqref{eq:cone}.

In view of the form of the filter function in~\eqref{eq:linearOptMu}, $\mu(\lambda)$ can be written as
\ba
\label{eq:muLinear}
\mu(\lambda)=\fv_a(\lambda)\sp{T}\tilde\xv_a+\fv_c(\lambda)\sp{T}\tilde\xv_c+\fv_s(\lambda)\sp{T}\tilde\xv_s,
\ea
where
\bse
\ba
\tilde\xv_a=\begin{bmatrix} \tilde a_1\\ \tilde a_2\\ \vdots \\ \tilde a_\Nf\end{bmatrix},\qquad 
\tilde\xv_c=\begin{bmatrix} \tilde c_1\\ \tilde c_2\\ \vdots \\ \tilde c_\Nf\end{bmatrix},\qquad 
\tilde\xv_s=\begin{bmatrix} \tilde s_1\\ \tilde s_2\\ \vdots \\ \tilde s_\Nf\end{bmatrix},
\ea
and
\ba
\label{eq:mucompV}
\fv_a(\lambda)=\begin{bmatrix} f_a(\lambda;\,T_{f,1})\\ f_a(\lambda;\,T_{f,2})\\ \vdots \\ f_a(\lambda;\,T_{f,\Nf})\end{bmatrix},\qquad
\fv_k(\lambda)=\begin{bmatrix} f_k(\lambda;\,\omega_1,T_{f,1})\\ f_k(\lambda;\,\omega_2,T_{f,2})\\ \vdots \\ f_k(\lambda;\,\omega_\Nf,T_{f,\Nf})\end{bmatrix},\quad \hbox{$k=c$ or $s$.}
\ea
\ese
The linear constraints in~\eqref{eq:constraintsT} can also be expressed in matrix form as
\ba
\label{eq:constraintsTM}
\Fmat_a\tilde\xv_a+\Fmat_c\tilde\xv_c+\Fmat_s\tilde\xv_s=\mathbf{1},\qquad
\Fmat_a\sp\prime\tilde\xv_a+\Fmat_c\sp\prime\tilde\xv_c+\Fmat_s\sp\prime\tilde\xv_s=\mathbf{0},
\ea
where the $\Nf\times\Nf$ matrices, $\Fmat_k$ and $\Fmat_k\sp\prime$, $k=a$, $c$ or $s$, are given by
\ba
\label{eq:constraintsVec}
\Fmat_k=\begin{bmatrix} \fv_k(\omega_1)\sp{T} \\ \fv_k(\omega_2)\sp{T} \\ \vdots \\ \fv_k(\omega_\Nf)\sp{T} \end{bmatrix},\qquad 
\Fmat_k\sp\prime=\begin{bmatrix} \fv_k\sp\prime(\omega_1)\sp{T} \\ \fv_k\sp\prime(\omega_2)\sp{T} \\ \vdots \\ \fv_k\sp\prime(\omega_\Nf)\sp{T} \end{bmatrix}.
\ea
Numerical experiments have shown it to be effective to use the $2\Nf$ constraints in~\eqref{eq:constraintsTM} to solve for $\tilde\xv_c$ and $\tilde\xv_s$ in terms of $\tilde\xv_a$.\footnote{Using $\tilde\xv_c$ or $\tilde\xv_s$ as free parameters can result in singular systems if the frequencies $\omega_m$ are integer multiples of one another.}  Accordingly, we write
\bse
\label{eq:yt}
\ba
\Gmat\tilde\yv=\zv-\Hmat\tilde\xv_a,
\ea
where
\ba
\Gmat=\begin{bmatrix} \Fmat_c & \Fmat_s \\ \Fmat_c\sp\prime & \Fmat_s\sp\prime \end{bmatrix},\qquad 
\tilde\yv=\begin{bmatrix} \tilde\xv_c  \\ \tilde\xv_s \end{bmatrix},\qquad 
\zv=\begin{bmatrix} \mathbf{1}  \\ \mathbf{0} \end{bmatrix},\qquad 
\Hmat=\begin{bmatrix} \Fmat_a  \\ \Fmat_a\sp\prime \end{bmatrix}. 
\ea
\ese
Assuming $\Gmat$ is nonsingular, we may solve for $\tilde\yv$ in terms of $\tilde\xv_a$, and then substitute the result into~\eqref{eq:muLinear} to give
\ba
\label{eq:muNice}
\mu(\lambda)=\gv_a(\lambda)\sp{T}\tilde\xv_a+\mu_0(\lambda),
\ea
where
\bse
\ba
\gv_a(\lambda)\sp{T}&=\fv_a(\lambda)\sp{T}-\bigl[\fv_c(\lambda)\sp{T}\;\;\fv_s(\lambda)\sp{T}\bigr]\Gmat\sp{-1}\Hmat, \\\mu_0(\lambda)&=\bigl[\fv_c(\lambda)\sp{T}\;\;\fv_s(\lambda)\sp{T}\bigr]\Gmat\sp{-1}\zv.
\ea
\ese

\newcommand{\lamMin}{{\lambda_{{\rm min}}}}
\newcommand{\lamMax}{{\lambda_{{\rm max}}}}

Using the result in~\eqref{eq:muNice}, with $\tilde\xv_c$ and $\tilde\xv_s$ given in terms of $\tilde\xv_a$ from~\eqref{eq:yt}, the optimization problem in~\eqref{eq:PoptT} and~\eqref{eq:constraintsT} becomes
\ba
\label{eq:xOpt}
\tilde\xv_{a,{\rm opt}}=\argmin_{\tilde\xv_a}\int_\Lambda \bigl(\gv_a(\lambda)\sp{T}\tilde\xv_a+\mu_0(\lambda)\bigr)\sp2\,d\lambda,
\ea
and then
\ba
\label{eq:yOpt}
\begin{bmatrix}\tilde\xv_{c,{\rm opt}} \\ \tilde\xv_{s,{\rm opt}}\end{bmatrix}=\tilde\yv_{{\rm opt}}=\Gmat\sp{-1}\bigl(\zv-\Hmat\tilde\xv_{a,{\rm opt}}\bigr).
\ea
In practice, the optimization problem is solved numerically by selecting a set of $\Nc_\lambda$ equally spaced points covering the interval $\Lambda$ with $\Nc_\lambda>\Nf$. 

{
\newcommand{\figw}{10.5cm}

\begin{figure}[htb] 
\label{fig:filterIllustration}
    \centering
    \begin{tikzpicture}

        \useasboundingbox (0,0.3) rectangle (\figw,3.8cm);

        \figByWidth{0}{0}{fig/filterGraph}{\figw}[0.][0.][0.][0.]

        \draw[<->,thick]
            (2.35,0.70) -- (2.95,0.70)
            node[midway,above,font=\scriptsize] {$\delta_1$};
        
        \draw[<->,thick]
            (4.00,.70) -- (4.78,.70)
            node[midway,above,font=\scriptsize] {$\delta_2$};

        \draw[<->,thick]
            (6.15,.70) -- (6.85,.70)
            node[midway,above,font=\scriptsize] {$\delta_3$};

        \draw[<->,thick]
            (6.15,1.50) -- (9.45,1.50)
            node[midway,above,font=\scriptsize ] {$c_{\delta}\,\ssf{\Nc_f}$$\delta_{\Nc_f}$};

        \draw[dashed,thick]
            (9.45,.70) -- (9.45,1.50); 

        \node[font=\scriptsize] at (1.15,0.75) {$0$};
        \node[font=\scriptsize] at (1.15,3.17) {$1$};   

    \end{tikzpicture}
    \caption{Illustration of the interval $\Lambda\in[\lambda_{{\rm min}},\lambda_{{\rm max}}]$ for the calculation of optimal filter parameters for a sample case with $\Nf=3$. The half-width of the peak in $\mu(\lambda)$ at $\omega_m$ is measured by $\delta_m$, $m=1,\ldots,\Nf$.}
\end{figure}
}

The interval $\Lambda = [\lamMin, \lamMax]$ is chosen with
\bse
\label{eq:LambdaRange}
\ba
\lamMin = 0,\qquad \lamMax = \omega_{\ssf\Nc_f} + c_{\ssf \delta}\Nc_f\;\delta_{\ssf \Nc_f},
\ea
where $c_{\ssf \delta}$ is a constant and $\delta_m$ denotes the half-width of the $m$-th peak of $|\mu|$, defined as
\ba
\delta_m \eqdef \f{\omega_m}{2N_{p,m}},
\ea
which is the distance from $\omega_m$ to the first zero in $\mu(\lambda)$ to the left of the $m$-th peak. The upper limit $\lamMax$ is thus chosen to extend beyond the last target frequency $\omega_{\ssf\Nc_f}$ by an amount proportional to $\delta_{\ssf \Nc_f}$, the half-width of the last peak about $\omega_{\ssf\Nc_f}$, see Figure~\ref{fig:filterIllustration}.
The choice in~\eqref{eq:LambdaRange} ensures that the optimization covers the full spectral region of interest including the tail beyond the last peak using $c_{\ssf \delta}=2$ (which our numerical experiments have shown to be a reasonable value).
The number of quadrature points $\Nc_\lambda$ is determined by specifying the number of points $\Nc_\delta$ per half-width $\delta_{\Nc_f}$, giving
\ba
\Nc_\lambda = \left\lceil \Nc_\delta \f{\lamMax}{\delta_{\Nc_f}} \right\rceil,
\ea
\ese
where $\lceil \cdot \rceil$ denotes the ceiling function, ensuring $\Nc_\lambda$ 
is a positive integer.  We set $\Nc_\delta=20$ for all calculations in this paper to adequately resolve the peaks of $\mu(\lambda)$.
A midpoint rule approximation of the integral in~\eqref{eq:xOpt} is used which leads to a standard least-squares problem to determine $\tilde\xv_{a,{\rm opt}}$, and then $\tilde\xv_{c,{\rm opt}}$ and $\tilde\xv_{s,{\rm opt}}$ are obtained from~\eqref{eq:yOpt}.  The steps to determine $\tilde\xv_{k,{\rm opt}}$, $k=a$, $c$ and $s$, are outlined in Algorithm~\ref{alg:SSWHV2}. See Appendix~\ref{sec:matlab_codes_to_determine_optimal_parameters} for Matlab codes implementing the algorithm.

\medskip

\renewcommand{\algFontSize}{\small}
\begin{algorithm}
\algFontSize 
\caption{Multi-Frequency WaveHoltz Filter Parameter Optimization.}
\begin{algorithmic}[1]

  \Function{\blue optParameters}{$\omega_m$, $m=1,2,\ldots,\Nf$}
  	\State Define $\fv_a(\lambda)$, $\fv_c(\lambda)$, $\fv_s(\lambda)$;  \Comment Given in~\eqref{eq:mucompV}
    \State Compute $\Fmat_k$, $\Fmat\sp\prime_k$, $k = a$, $c$ or $s$; \Comment Given in~\eqref{eq:constraintsVec}
    \State Assemble $\Gmat, \; \Hmat, \; \zv$; \Comment Given in~\eqref{eq:yt}
    \State $\Lambda = [\lamMin, \lamMax]$; \Comment Select the spectral interval, see~\eqref{eq:LambdaRange}  
    \State Select quadrature points $\{\lambda_k\}_{k=1}^{\Nc_\lambda}$ covering $\Lambda$; \Comment $\Nc_\lambda$ equally spaced points
    \State Evaluate $\mathbf{g}_a(\lambda_k)^T = \mathbf{f}_a(\lambda_k)^T - [\mathbf{f}_c(\lambda_k)^T \;\; \mathbf{f}_s(\lambda_k)^T] \Gmat^{-1} \Hmat$;  \Comment for $k=1,\ldots,\Nc_\lambda$
    \State Evaluate $\mu_0(\lambda_k) = [\mathbf{f}_c(\lambda_k)^T \;\; \mathbf{f}_s(\lambda_k)^T] \Gmat^{-1}\zv$; \Comment for $k=1,\ldots,\Nc_\lambda$
    \State Set $\Mmat = [\mathbf{g}_a(\lambda_k)^T]_{k=1}^{\Nc_\lambda}, \quad \mv_0 = -[\mu_0(\lambda_k)]_{k=1}^{\Nc_\lambda}$; \Comment Assemble the least-squares system
    \State $\tilde{\xv}_{a,\text{opt}} = \Mmat\,\backslash\,\mathbf{m}_0$  \Comment  Solve the least-squares problem.
    \State $ \begin{bmatrix} \tilde{\xv}_{c,\text{opt}} \\ \tilde{\xv}_{s,\text{opt}} \end{bmatrix} = \tilde{\yv}_{\text{opt}} = \Gmat\sp{-1}\left(\mathbf{z} - \Hmat\tilde{\xv}_{a,\text{opt}}\right)$ \Comment Determine $\tilde{\xv}_{c,\text{opt}}$ and $\tilde{\xv}_{s,\text{opt}}$
 \EndFunction
\end{algorithmic}
\label{alg:SSWHV2}
\end{algorithm}

As a first example, consider again the case $\Nf=2$ with $\omega_1=5$ and $\omega_2=9$. The least-squares solution yields the parameters
{
\sisetup{group-digits=false}
\begin{table}[H]\footnotesize
	\begin{center}
	\setlength{\tabcolsep}{4.5pt}   
		\begin{tabular}{|c|c|c|S[table-format=-1.6]S[table-format=-1.6]S[table-format=-1.6]|S[table-format=-1.6]S[table-format=-1.6]S[table-format=-1.6]|}
			\hline
			\multirow{2}{*}{$N_{p,1}$} & \multirow{2}{*}{$N_{p,2}$} & \multirow{2}{*}{$\lamMax$}
			& \multicolumn{3}{c}{$\tilde\pv_1 = (\tilde a_1, \tilde c_1, \tilde s_1)$}
			& \multicolumn{3}{c|}{$\tilde\pv_2 = (\tilde a_2, \tilde c_2, \tilde s_2)$} \\
			\cline{4-9}
			&  &  &{$\tilde a_1$} & {$\tilde c_1$} & {$\tilde s_1$} & {$\tilde a_2$} & {$\tilde c_2$} & {$\tilde s_2$} \\
			\hline
			1 & 1 & 27   & 0.9054    & 0.03864  & 1.655     & -1.107  & 1.324  & 0.7038  \\
			2 & 3 & 15   & -0.4327  & 0.9911    & -0.007299 & 0.3830  & 1.051  & -0.02323 \\
			3 & 5 & 12.6 & -0.07091  & 0.9684   & 0.05086   & 0.07765 & 0.9583 & -0.1045 \\
			\hline
		\end{tabular}
	\end{center}
	\caption{Optimal parameter values for $\omega_1 = 5$, $\omega_2 = 9$. Here $N_p = N_{p,1}$.}
	\label{tab:params}
\end{table}
}

Figure~\ref{fig:MFWHfilterFunction2b} plots the behavior of $\mu$ for $N_p=1$, $2$ and~$3$ for this two-frequency case with $\omega_1=5$ and $\omega_2=9$ using the parameters given in Table~\ref{tab:params} obtained using Algorithm~\ref{alg:SSWHV2}. Here $\beta_m$ evaluated at $\omega_m$ is not necessarily equal to one; rather, the parameters $\tilde\pv_m$ are chosen such that $\mu(\omega_m)=1$ and $\lambda=\omega_m$ is a local extremum, and so $\tilde\pv_m$ may differ for each $N_p$.
The three plots in the figure from left to right show the behavior of $\mu$ for $N_p = 1$, $2$ and $3$, respectively.  As before, we note that the peaks sharpen as $N_p$ increases, but unlike the corresponding plots in Figure~\ref{fig:MFWHfilterFunction2} which use the traditional filter parameters, we now see that there are no intervals where $|\mu| > 1$.  Thus, the FPI converges for each choice of $N_p$. This also avoids the possibility of  spurious resonances, ensuring that the matrix forming the Krylov basis in GMRES is nonsingular and that GMRES converges to the correct solution.

{
\newcommand{\figSize}{5.3cm}
\begin{figure}[htp]
\begin{center}
\resizebox{\textwidth}{!}{
\begin{tikzpicture}[scale=1]
  \useasboundingbox (0,.35) rectangle (3*\figSize,4.5);  

   \begin{scope}[yshift=0cm]
    \figByWidth{0*\figSize}{0}{fig/optmfwhOmega5Omega9Np1BddddFixPointMuFunction}{\figSize}[.0][0.0][0.][0.]
    \figByWidth{1*\figSize}{0}{fig/optmfwhOmega5Omega9Np2BddddFixPointMuFunction}{\figSize}[.0][0.0][0.][0.]
    \figByWidth{2*\figSize}{0}{fig/optmfwhOmega5Omega9Np3BddddFixPointMuFunction}{\figSize}[.0][0.0][0.][0.] 
       
  \end{scope}  
\end{tikzpicture}
}
\end{center}
\caption{Multi-frequency WaveHoltz filter function for two frequencies,  $\omega_1=5$ and
$\omega_2=9$. The solution is integrated over $\Np$ periods of the longest period $T_1=2\pi/\omega_1$.
    } 
\label{fig:MFWHfilterFunction2b}
\end{figure}
}

{
	\newcommand{\figSize}{5.3cm}
	\begin{figure}[htp]
		\begin{center}
			\resizebox{\textwidth}{!}{
				\begin{tikzpicture}[scale=1]
					\useasboundingbox (0,.35) rectangle (3*\figSize,4.5);  

					\begin{scope}[yshift=0cm]
						\figByWidth{0*\figSize}{0}{fig/OriginaloptmfwhOmega5Omega9Np1BddddFixPointMuFunction}{\figSize}[.0][0.0][0.][0.]
						\figByWidth{1*\figSize}{0}{fig/OriginaloptmfwhOmega5Omega9Np2BddddFixPointMuFunction}{\figSize}[.0][0.0][0.][0.]
						\figByWidth{2*\figSize}{0}{fig/OriginaloptmfwhOmega5Omega9Np3BddddFixPointMuFunction}{\figSize}[.0][0.0][0.][0.] 
						
					\end{scope}  
				\end{tikzpicture}
			}
		\end{center}
		\caption{Multi-frequency WaveHoltz filter functions for two frequencies, $\omega_1 = 5$ and $\omega_2 = 9$, comparing the original filter $\mu_{org}$ and the optimized filter $\mu_{opt}$. The solution is integrated over $N_p$ periods of the longest period $T_1 = 2\pi/\omega_1$.
		}
		\label{fig:MFWHfilterFunction2c}
	\end{figure}
}

Figure~\ref{fig:MFWHfilterFunction2c} compares the optimized filter $\mu_{{\rm opt}}$ with the original (traditional) filter $\mu_{{\rm org}}$ for the same frequencies $\omega_1 = 5$ and $\omega_2 = 9$, for $N_p = 1, 2$, and 
$3$. For $N_p = 1$ (left), the two filters differ noticeably: $\mu_{{\rm org}}$ exhibits larger intervals where $|\mu|>1$, whereas $\mu_{{\rm opt}}$ achieves a sharper peak near the target frequencies. As $N_p$ increases to $2$ and~$3$ (middle and right), the two filters become almost indistinguishable near the target frequencies, with the optimization offering diminishing improvement as $N_p$ grows. This is consistent with the fact that for sufficiently large $N_p$ the peaks of $\beta_1$ and $\beta_2$ are already well separated, leaving little room for the least-squares optimization to further improve the filter, except near $\lambda = 0$.

\newcommand{\omegaAdjust}{\tilde{\omega}}
\section{Discretizing the wave equation and filter} \label{sec:discretizingTheWaveEquation}

The MFWH algorithm can be implemented with any number of discrete methods such as those based on finite element, 
finite difference and finite volume approximations. 
Here we consider using finite difference methods which can be extended to complex geometries using overset grids~\cite{overHoltzPartOne}.
Consider solving the IBVP for the wave equation given by~\eqref{eq:waveEqnMultipleForcings} on a 
single mapped (structured) grid.
Let $\WFunc_{\jv}^n \approx \wFunc(\xv_{\jv},t^n)$ denote the approximate solution at
grid point $\xv_{\jv}$ and time $t^n=n\dt$, where $\dt$ is the time-step. 
Here $\jv=[j_1,j_2,j_3]$ is a multi-index with $j_l=0,1,2\ldots,N_l$ denoting the range of indices including the boundaries and interior points, 
where $N_l$ is the number of grid cells in direction $l$. Additional ghost points are added as needed.
 We consider both explicit and implicit methods in time. 
While high-order accurate methods in both space and time are available (see~\cite{lcbc2022} for example), 
we use only second-order accurate schemes in time since we can correct for time discretization errors in the
 WaveHoltz algorithm as discussed below.
The spatial approximations, on the other hand, are $p^{\rm th}$ order accurate, where $p=2$ and $4$ for the purposes of this paper (although higher-order accurate discretizations are possible).

\subsection{Explicit and implicit time-stepping} \label{sec:timeStepping}

The discrete approximation to the second time derivative in the wave equation uses a standard second-order accurate central difference, which, in terms of the
forward and backward divided difference operators, $\Dpt \WFunc_\jv^n \eqdef  (\WFunc_\jv^{n+1}-\WFunc_\jv^{n})/\dt$ and $\Dmt \WFunc_\jv^n  \eqdef  (\WFunc_\jv^{n}-\WFunc_\jv^{n-1})/\dt$, respectively, is 
$
    \Dpt\Dmt \WFunc_\jv^n  = ( \WFunc_\jv^{n+1} -2  \WFunc_\jv^n +  \WFunc_\jv^{n-1} )/(\dt^2). 
$
The explicit or implicit time-stepping scheme for the wave equation with modified frequencies
 $\omegaAdjust_m$ (chosen to adjust for time-discretization errors)
takes the form 
\bse
\label{eq:WaveScheme}
\bat
   & \Dpt\Dmt \WFunc_\jv^n =  L_{ph} \Big( \alphaI \WFunc_\jv^{n+1} + \betaI \WFunc_\jv^n + \alphaI \WFunc_\jv^{n-1}  \Big) - F_\jv^n , 
         \quad && \jv \in \Omega_h, ~ n=0,1,\ldots  ,   \label{eq:WaveSchemeA}  
\eat         
\bat
   & F_\jv^n \eqdef  \sum_{m=1}^{\Nf} f^{(m)}(\xv_\jv) \, \cos( \omegaAdjust_m t^n)  \,(\betaI+2\alphaI\cos(\omegaAdjust_m\dt)),  \label{eq:WaveSchemeForcing}\\
   & \WFunc_\jv^0 = \sum_{m=1}^\Nf V_\jv^{(m,k)},                                                        \qquad&& \jv \in \bar{\Omega}_h  , \\
   & \Dzt \WFunc_\jv^0 = 0,                                                             \qquad&& \jv \in \bar{\Omega}_h  ,  \label{eq:WaveSchemeC} \\
   & \Bc_{ph} \, \WFunc_\jv^n = G_\jv^n \eqdef \sum_{m=1}^\Nf g^{(m)}(\xv_\jv) \, \cos( \omegaAdjust_m t^n)  , \qquad&&  \jv \in \p\Omega_h , ~ n=1,2,\ldots ~. \label{eq:WaveSchemeBC}
\eat
\ese
Here $L_{ph}$ is a $p$-th order accurate approximation to $L$, $\Bc_{ph}$ denotes the discretized boundary condition operator, 
$\Omega_h$ denotes the set of grid points, $\jv$, where the interior equation is applied,
$\bar{\Omega}_h$ denotes the set of all grid points, and $\p\Omega_h$ denotes the set of points
where the boundary conditions are applied. 
The constants $(\alphaI,\betaI)$ in~\cref{eq:WaveSchemeA} satisfy $\betaI=1-2\alphaI$ for second-order accuracy in time. 
The explicit time-stepping scheme takes $\alphaI=0$ and $\betaI=1$.  There are many choices for implicit schemes with $\alphaI>0$, and we use $\alphaI=1/2$ and $\betaI=0$ corresponding to the trapezoidal method which is unconditionally stable (see~\cite{wimp2025} for further discussion).
The modified frequencies $\omegaAdjust_m$ appearing in~\cref{eq:WaveSchemeForcing} and~\cref{eq:WaveSchemeBC} are 
\ba
    \omegaAdjust_m \eqdef \f{1}{\dt} \cos^{-1}\Big( \f{1- (\betaI/2) (\omega_m\dt)^2}{1 + \alphaI(\omega_m\dt)^2 } \Big) ,
    \label{eq:omegaImplicit}
\ea
with corresponding modified period $\Ttilde_m=(2\pi)/\omegaAdjust_m$.
The initial condition~\eqref{eq:WaveSchemeC} can be combined with~\eqref{eq:WaveSchemeA}
when $n=0$ to eliminate $\WFunc_\jv^{-1}$, and thus arrive at an implicit update for the first time-step~$\WFunc_\jv^{1}$, 
\ba
   \Big[ I - \alphaI \dt^2 L_{ph} \Big] \WFunc_\jv^{1} & = 
       \WFunc_\jv^0  + \f{\betaI}{2} \dt^2 L_{ph} \WFunc_\jv^0 
        - \half \dt^2  F_\jv^0 .
       \label{eq:implicitFirstStep}
\ea
Fortuitously the implicit matrix implied by the left-hand side of~\cref{eq:implicitFirstStep} is the same one implied by~\cref{eq:WaveSchemeA}, and thus the same implicit solver used for the interior scheme 
can be used on the first step but with an altered right-hand side.
Using the adjusted frequencies defined in~\cref{eq:omegaImplicit}, the iterates $V_\jv^{(m,k)}$ converge, as $k\rightarrow \infty$, to 
the discrete Helmholtz solutions $U_\jv^{(m)}$ satisfying
\bse
\label{eq:discreteHelmholtz} 
\bat
   & L_{ph} U_\jv^{(m)} + \omega_m^2 U_\jv^{(m)} = f^{(m)}(\xv_\jv),  \qquad && \jv \in \Omega_h,       \\
   & \Bc_{ph} \, U_\jv^{(m)} = g^{(m)}(\xv_\jv)  ,                         \qquad&&  \jv \in \p\Omega_h ,
\eat
\ese
with frequencies $\omega_m$ agreeing with the exact problem in~\eqref{eq:MultipleHelmholtz}.
\subsection{Discrete time filter}

The time filters in~\cref{eq:pDef} are approximated with numerical integration using the modified frequencies and periods, 
\ba
\label{eq:pDefQ}
   & \rhs_{m,\jv} = \f{2}{\Ttilde_{f,m}} \sum_{n=0}^{N_{t,m}} \sigma_{n,m}\, \filter(\omegaAdjust_mt;\params) \, \WFunc_\jv^n  ,
\ea
where $\sigma_{n,m}$ are quadrature weights and $N_{t,m}$ denotes the number of time-steps involved in the quadrature.  The parameters, $\params$, in~\eqref{eq:pDefQ} are taken to be the traditional values, $\pv_0$ (for all $m$), or the least-squares parameters, $\paramsT_m$.
For the longest period, a trapezoidal rule approximation is used. 
The quadrature weights for other frequencies can be easily derived as follows. Consider developing a numerical integration rule in time for some function~$\zeta(t)$,
\ba\label{eq:DefNI}
   \int_{0}^{\Ttilde_{f,m}} \zeta(t)  dt \approx \sum_{n=0}^{N_t} \sigma_{n,m} \, \zeta(t^n) ,
\ea
where $N_t$ denotes the total number of time-steps.\footnote{It is convenient in a code to keep all terms in the sum with some $\sigma_{n,m}$ set to zero.} Suppose the final time $\Ttilde_{f,m}$
lies in the time interval $[t^{q},t^{q+1}]$ for some integer 
$q<N_t$.  
At second-order accuracy a composite trapezoidal type rule can be used where 
a linear polynomial, using $\zeta(t^q)$ and $\zeta(t^{q+1})$,  is employed to approximate $\zeta(t)$ for the interval $[t^{q},t^{q+1}]$. This polynomial is integrated over the sub-interval $[t^q,\Ttilde_{f,m}]$ to give adjusted weights 
at $\sigma_{q,m}$ and $\sigma_{q+1,m}$. Weights $\sigma_{n,m}$ are set to zero for $n>q+1$.

The discrete beta functions using these integration rules are 
 \ba
    \beta_{d,m}(\lambda;\omegaAdjust_m,\Ttilde_{f,m},\params) \eqdef \f{2}{\Ttilde_{f,m}} \sum_{n=0}^{N_{t,m}} \sigma_{n,m}  \,\filter(\omegaAdjust_m t^n;\params) \cos(\lambda t^n) . 
    \label{eq:betadDefA}
 \ea
 Note that if the parameters $\pv_m$ are chosen to follow the original single-frequency filter, then the nominal choice is $\pv_0=[-1/4,1,0]$.
In the discrete case, this choice is modified to $[a_m,1,0]$, where $a_m$ is specified by setting the derivative of $\beta_{d,m}$ with respect to $\lambda$ equal to zero when $\lambda=\omegaAdjust_m$. The result is
\ba
\label{eq:origParamAdjust}
   a_m = - \f{\tan(\omegaAdjust_m\dt/2)}{2 \tan(\omegaAdjust_m\dt)}.  
\ea
The time-discrete version of $\mu$ in~\cref{eq:MFWHmu} is denoted by $\mu_d$, and is defined in terms of $\beta_{d,m}$ for the original parameters as
\ba
  \mu_d(\lambda) \eqdef \sum_{m=1}^{\Nf}\weight_m \, \beta_{d,m}(\lambda) . \label{eq:MFWHmud}
\ea
Note that no adjustments are made to the parameters $\paramsT_m$ resulting from the least-squares optimization.  For this case, we have
\ba
  \mu_d(\lambda) \eqdef \sum_{m=1}^{\Nf} \tilde{\beta}_{d,m}(\lambda) , \label{eq:MFWHmudN}
\ea
where
 \ba
    \tilde{\beta}_{d,m}(\lambda;\omegaAdjust_m,\Ttilde_{f,m},\paramsT_m) \eqdef \f{2}{\Ttilde_{f,m}} \sum_{n=0}^{N_{t,m}} \sigma_{n,m}  \,\filter(\omegaAdjust_m t^n;\paramsT_m) \cos(\lambda t^n) . 
    \label{eq:betadDefAT}
 \ea

\subsection{Time-step determination}

For explicit time-stepping the maximum stable time-step is determined by a von Neumann analysis and then a safety factor ${\rm CFL}=0.9$ is applied.
For implicit time-stepping the trapezoidal or full-weighting schemes have no time-step restriction. The convergence of the WaveHoltz iteration, however,
does require there to be at least five time-steps per period~\cite{overHoltzPartOne}. In practice we usually use $10$ time-steps over the 
smallest period $T_{f,\Nf}$ to estimate $\dt$ and this leads to the total number of time-steps used to reach the overall final time $T_{f,1}$.
Note that the adjusted frequencies and periods depend on $\dt$, which itself depends on the final period; the adjusted final time and time-step 
to satisfy these constraints can be found following the recipe in~\cite{overHoltzPartOne}.

\section{Acceleration by Krylov methods} \label{sec:algorithm_gmres}

The discrete form of the WaveHoltz mapping~\eqref{eq:waveHoltzMapping} can be written as 
\ba
    \Vv_h^{(k+1)} = \Wc_h \,\Vv_h^{(k)} = S_h \, \Vv_h^{(k)} + \bv_h , \label{eq:discreteMFWHmapping}
\ea
where $S_h$ is a matrix, and where $\Vv_h^{(k)}$ is the vector containing the grid functions $V_\jv^{(m,k)}$, $m=1,2,\ldots,\Nf$. 
The vector $\bv_h$, which contains information from the forcings, 
is the result of taking one step of the MFWH algorithm with zero initial conditions, $\bv_h = \Wc_h \zerov$.
At convergence~\cref{eq:discreteMFWHmapping} is equivalent to solving the linear system
\ba
\label{eq:krylovSystem}
    A_h \Vv_h = \bv_h, \qquad A_h \eqdef I - S_h. 
\ea
As shown in Section~\ref{sec:numerical_confirmation_of_the_convergence_analysis}, the rate of convergence is typically much faster if a Krylov space method is used to solve the linear system in~\eqref{eq:krylovSystem} directly rather than using the fixed-point iteration implied by~\eqref{eq:discreteMFWHmapping}.  Since $A_h$ is not symmetric for the multi-frequency case, in general, we use GMRES, which is shown to be an effective choice.
An additional advantage of using a Krylov solver, such as GMRES, 
is that the existence of discrete eigenvalues $\lambda_{h,\nu}$ in intervals where $| \mu_d | > 1$ does not prevent convergence.

Note, however, if an interval (or intervals) of $\lambda$ exist where $\mu_d(\lambda)>1$, then an eigenvalue $\lambda_{h,\nu}\ne\omega_m$, $m=1,2,\ldots,\Nf$, can occur with $\mu_d(\lambda_{h,\nu})=1$ 
or nearly equal to one (i.e.~a \textsl{spurious resonance} or \textsl{near spurious resonance} can occur).
As an example, see the left graph in Figure \ref{fig:MFWHfilterFunction2} where there are two spurious resonances for $\lambda$ near $7$ and $11$.
In this case there could be an eigenvector of $S_h$ with eigenvalue equal to one (or close to 1) and thus the corresponding eigenvector of $A_h$ has eigenvalue equal to zero (or close to zero). 
The matrix $A_h$ would therefore be singular (or nearly singular) and the Krylov solver would likely fail to give the correct answer. 
This failure of the GMRES algorithm can be detected a~posteriori by checking that the residual in the original discretized Helmholtz equations is small.

The issue of a singular Krylov matrix is discussed further in Section~\ref{sec:numerical_confirmation_of_the_convergence_analysis}.  It is shown that a filter function for which  spurious resonances occur leads to the MFWH algorithm with GMRES converging to a non-solution.  If, on the other hand, filter parameters from the least-squares problem are used, then the filter function for the same problem configuration has no spurious resonances and the MFWH algorithm with GMRES converges to the true Helmholtz solutions without difficulty.

\newcommand{\Nit}{N_{\rm it}}
\section{Numerical results and $O(N)$ scaling} \label{sec:numerical_confirmation_of_the_convergence_analysis}

Theorem~\ref{thm:MFWHmu} provides the condition for the convergence of the MFWH~fixed-point iteration, and in this
section some numerical results are given to confirm the theory.
The discrete asymptotic convergence rate, $\ACR$, is computed as 
\ba
    \ACR \eqdef \max_\nu |\mu_d(\lambdaTilde_{h,\nu})|,
\ea
where $\mu_d$ (given in~\cref{eq:MFWHmud} or~\cref{eq:MFWHmudN}) is the discrete version of the WaveHoltz filter function $\mu$ 
and $\lambdaTilde_{h,\nu}$ are eigenvalues of the discrete version of the eigenvalue problem~\cref{eq:eigBVP} that have been adjusted to account
for time discretization errors using 
\ba
  \lambdaTilde(\lambda) \eqdef \f{1}{\dt} \cos^{-1}\Big( \f{1- (\betaI/2) (\lambda\dt)^2}{1 + \alphaI(\lambda\dt)^2 } \Big) .
\ea
The average convergence rate for the MFWH iteration over a total of $\Nit$ iterations is defined as
\ba
    \CR \eqdef \left[\frac{ r^{(\Nit)} }{ r^{(1)} }\right]^{1/\Nit} ,
\ea
where the (combined) \textsl{residual} in the $k\sp{{\rm th}}$ iterate $r^{(k)}$ is computed from
\ba 
    (r^{(k)})^2            & \eqdef \f{1}{\Nf N_a} \sum_{m=1}^{\Nf} \, \sum_{\jv\in\Omega_{h}}  \big( V_\jv^{m,k} - V_\jv^{m,k-1} \big)^2,  \label{eq:residual}
\ea
where $N_a$ is the number of grid points in $\Omega_h$. 
In order to more fairly compare the convergence rates for different values of $\Np$ we define
the effective  convergence rate as
\ba
  \ECR = \CR^{1/\Np}, 
\ea
which takes into account the extra work as $\Np$ increases. Thus, for example, if $\CR=0.5$ for $\Np=1$ and $\CR=0.25 = (0.5)^2$ for $\Np=2$,
then $\ECR=0.5$ in both cases.

The numerical examples use Gaussians for the source terms in~\eqref{eq:MultipleHelmholtz} given by
\ba
   f^{(m)}(\xv) = a_{g,m} \exp\bigl( - b_{g,m}^2 \, \| \xv-\xv_{0,m}\|^2 \bigr),
   \label{eq:gaussianSource}
\ea
where $a_{g,m}$ is the amplitude, $\xv_{0,m}=[x_{0,m},y_{0,m}]$ denotes the center of the Gaussian, and $b_{g,m}^{-1}$ determines the approximate width of the Gaussian. The harmonic time dependence of the terms in the composite source $F(\xv, t)$ in~\eqref{eq:waveEqnCompositeForcings} use the adjusted frequencies $\omegaAdjust_m$ as discussed in Section~\ref{sec:timeStepping}.
Dirichlet boundary conditions with $g^{(m)}(\xv)=0$ are used for all example problems.
The wave speed is taken to be one, i.e.~$c=1$, for all computations.

{
  \newcommand{\figw}{6cm}
  \newcommand{\figh}{5.2cm}
\begin{figure}[htb]
\begin{center}
\begin{tikzpicture}
  \useasboundingbox (0,.35) rectangle (13,.95*\figh);  

  \begin{scope}[yshift=0*\figh]
    \figByWidth{        0}{0}{fig/square128FD24TSENf3Np2}{\figw}[0.][0.][0.][0.]
    \figByWidth{1.1*\figw}{0}{fig/square128FD24TSENf3Np2FixPointMuFunctionDiscreteBeta}{\figw}[0.][0.][0.][0.]
  \end{scope}
\end{tikzpicture}
\end{center}
\caption{Multi-frequency algorithm for $\omega_m=\{5.1,10.1,15.1\}$ with $N_{p,m}=\{2,3,5\}$ and traditional filter parameters.  Explicit time-stepping and fourth-order accuracy on a square128 grid.
   Left: Behavior of the residuals for the FPI and GMRES iterations. Right: MFWH fixed-point convergence function $|\mu_d(\lambdaTilde(\lambda))|$.
   The black vertical lines on the right graph mark the locations of $\omega_m$,
   while dashed lines denote the adjusted frequencies $\omegaAdjust_m$ (here almost identical).
   The {\red red x's} on the right graph mark locations of the eigenvalues while the
   value of $|\mu_d|$ at the {\green green circle} marks the ACR. }
\label{fig:multiFreq3SquareNp2Exp}
\end{figure}
}

{
  \newcommand{\figw}{6cm}
  \newcommand{\figh}{5.2cm}
\begin{figure}[htb]
\begin{center}
\begin{tikzpicture}
  \useasboundingbox (0,.35) rectangle (13,.95*\figh);  

  \begin{scope}[yshift=0*\figh]
    \figByWidth{        0}{0}{fig/square128FD24TSINf3Np2}{\figw}[0.][0.][0.][0.]
    \figByWidth{1.1*\figw}{0}{fig/square128FD24TSINf3Np2FixPointMuFunctionDiscreteBeta}{\figw}[0.][0.][0.][0.]
  \end{scope}
\end{tikzpicture}
\end{center}
\caption{Multi-frequency algorithm for $\omega_m=\{5.1,10.1,15.1\}$ with $N_{p,m}=\{2,3,5\}$ and traditional filter parameters.  Implicit time-stepping and fourth-order accuracy on a square128 grid.
   Left: Behavior of the residuals for the FPI and GMRES iterations. Right: MFWH fixed-point convergence function $|\mu_d(\lambdaTilde(\lambda))|$.
   The black vertical lines on the right graph mark the locations of $\omega_m$,
   while dashed lines denote the adjusted frequencies $\omegaAdjust_m$.
}
\label{fig:multiFreq3SquareNp2Imp}
\end{figure}
}

Figure~\ref{fig:multiFreq3SquareNp2Exp} shows results from the MFWH\ algorithm for solving
three different Helmholtz problems with frequencies $\omega_m = \{ 5.1, 10.1, 15.1 \}$ on a two-dimensional unit square domain.
The domain is covered by a Cartesean grid with $N_1=N_2=128$ grid cells, and the filter parameters are taken as the traditional choice, $\pv=\pv_0$, with adjustments following~\eqref{eq:origParamAdjust}.
The Gaussian source parameters in~\cref{eq:gaussianSource} are taken as
\ba
\label{eq:Nf3parameters}
\left\{ \begin{array}{l}
   a_{g,m} =  \{ 25, 100, 225 \}, \smallskip\\
   b_{g,m} =  \{ 15, 15, 15 \},  \smallskip \\
  \xv_{0,m} = \{ [0.6,0.45], [0.4,0.5], [0.55,0.5]   \} .
\end{array}\right.
\ea
Explicit time-stepping is utilized with fourth-order accuracy in space.
(The convergence results using second-order accuracy are very similar.) The converged MFWH solutions match direct solutions
of the three discrete Helmholtz problems to a relative error of about $10^{-12}$, see the definition of $e_h\sp{(m)}$ defined below in~\eqref{eq:errorAndResid}.
The left graph shows results for $N_p=2$ (i.e.~integrating over two periods of $T_1$). 
The right plot in the figure shows 
a graph of $|\mu_d(\lambdaTilde(\lambda))|$ versus~$\lambda$ and also plots red x's at
$|\mu_d(\lambdaTilde_{h,\nu})|$, where $\lambda_{h,\nu}$ are the adjusted eigenvalues of the Laplacian on the square.
A green circle marks the value of $|\mu_d|$ at the worst case eigenvalue and this value determines the ACR.
The convergence of the FPI is seen to match the theory. The iterations that employ GMRES are seen to converge rapidly.  Even though there are intervals for which $|\mu_d|>1$, there are no spurious resonances for this problem.

The graphs of Figure~\ref{fig:multiFreq3SquareNp2Imp} show corresponding results for implicit time-stepping using $10$ time-steps per-period for the smallest period.
The convergence results for the case of implicit time-stepping are very similar to those using explicit time-stepping.

{
	\newcommand{\figw}{6cm}
	\newcommand{\figh}{5.2cm}
	\begin{figure}[htb]
		\begin{center}
			\begin{tikzpicture}
				\useasboundingbox (0,.35) rectangle (13,.95*\figh);  

				\begin{scope}[yshift=0*\figh]
					\figByWidth{        0}{0}{fig/Sq128originalFD4TSINf3Np1}{\figw}[0.][0.][0.][0.]
					\figByWidth{1.1*\figw}{0}{fig/Sq128originalFD4TSINf3Np1FixPointMuFunctionDiscreteBeta}{\figw}[0.][0.][0.][0.]
				\end{scope}
			\end{tikzpicture}
		\end{center}
		\caption{Multi-frequency algorithm for $\omega_m=\{3,6,12\}$ with $N_{p,m}=\{1,1,3\}$ and traditional filter parameters.  Implicit time-stepping and fourth-order accuracy on a square128 grid.
   Left: Behavior of the residuals for the FPI and GMRES iterations. Right: MFWH fixed-point convergence function $|\mu_d(\lambdaTilde(\lambda))|$.
   The black vertical lines on the right graph mark the locations of $\omega_m$,
   while dashed lines denote the adjusted frequencies $\omegaAdjust_m$.
		}
		\label{fig:multiFreq3SquareNp1ImpEx3}
	\end{figure}
}

{
	\newcommand{\figw}{6cm}
	\newcommand{\figh}{5.2cm}
	\begin{figure}[htb]
		\begin{center}
			\begin{tikzpicture}
				\useasboundingbox (0,.35) rectangle (13,.95*\figh);  

				\begin{scope}[yshift=0*\figh]
					\figByWidth{        0}{0}{fig/Sq128leastsquareFD24TSINf3Np1}{\figw}[0.][0.][0.][0.]
					\figByWidth{1*\figw}{0}{fig/Sq128leastsquareFD24TSINf3Np1FixPointMuFunctionDiscreteBeta}{\figw}[0.][0.][0.][0.]
				\end{scope}
			\end{tikzpicture}
		\end{center}
		\caption{Multi-frequency algorithm for $\omega_m=\{3,6,12\}$ with $N_{p,m}=\{1,1,3\}$ and optimized filter parameters.  Implicit time-stepping and fourth-order accuracy on a square128 grid.
   Left: Behavior of the residuals for the FPI and GMRES iterations. Right: MFWH fixed-point convergence function and note that $|\mu_d(\lambdaTilde(\lambda))|\le1$.
   The black vertical lines on the right graph mark the locations of $\omega_m$,
   while dashed lines denote the adjusted frequencies $\omegaAdjust_m$.
		}
		\label{fig:multiFreq3SquareNp1ImpEx4}
	\end{figure}
}

{
	\newcommand{\figw}{6cm}
	\newcommand{\figh}{5.2cm}
	\begin{figure}[htb]
		\begin{center}
			\begin{tikzpicture}
				\useasboundingbox (0,.35) rectangle (13,.95*\figh);  

				\begin{scope}[yshift=0*\figh]
					\figByWidth{        0}{0}{fig/Sq128leastsquareFD24TSINf3Np2}{\figw}[0.][0.][0.][0.]
					\figByWidth{1.1*\figw}{0}{fig/Sq128leastsquareFD24TSINf3Np2FixPointMuFunctionDiscreteBeta}{\figw}[0.][0.][0.][0.]
				\end{scope}
			\end{tikzpicture}
		\end{center}
		\caption{Multi-frequency algorithm for $\omega_m=\{3,6,12\}$ with $N_{p,m}=\{2,3,7\}$ and optimized filter parameters.  Implicit time-stepping and fourth-order accuracy on a square128 grid.
   Left: Behavior of the residuals for the FPI and GMRES iterations. Right: MFWH fixed-point convergence function and note that $|\mu_d(\lambdaTilde(\lambda))|\le1$.
   The black vertical lines on the right graph mark the locations of $\omega_m$,
   while dashed lines denote the adjusted frequencies $\omegaAdjust_m$.
		}
		\label{fig:multiFreq3SquareNp1ImpEx5}
	\end{figure}
}

Figure~\ref{fig:multiFreq3SquareNp1ImpEx3} shows results for a case with frequencies $\omega_m = \{3, 6, 12\}$ and using the Gaussian source parameters given in~\eqref{eq:Nf3parameters}.  For this case, the MFWH\ algorithm is used with $\Np=1$ and filter parameters $\pv_0$ (with adjustments for time discretization errors as before). 
We observe that the FPI diverges (at a rate predicted by the theory) while GMRES converges quite rapidly. 
Figure~\ref{fig:multiFreq3SquareNp1ImpEx4} show results for the same case, but with filter parameters
\ba
\left\{ \begin{array}{l}
\tilde\pv_1 = (0.5111, 0.5737, 1.374),\smallskip\\
\tilde\pv_2 = (-0.9718, 0.9409, 0.5267),  \smallskip\\ 
\tilde\pv_3 = (0.3236, 1.184, 0.3190),
\end{array}\right.
\ea
obtained via least-squares optimization using $\Lambda = [0, 22.23]$, $\Nc_\lambda = 240$, $N_t = 40$ and $T_f = 2.105$. The optimal filter parameters obtained here, and in subsequent results, are computed with the Matlab code in Listing~\ref{code:optparams}. The optimized filter successfully reduces the dominant eigenvalue below 1, and  therefore the FPI now converges (at the rate predicted by the theory). GMRES also converges rapidly in this case.

In order to show the effect of increasing $\Np$, we use the MFWH\ algorithm for the case $\omega_m = \{3, 6, 12\}$ as before, but now using $\Np=2$ with parameters
\ba
\left\{ \begin{array}{l}
\tilde\pv_1 = (0.5111, 0.5737, 1.374),\smallskip\\
\tilde\pv_2 = (-0.9718, 0.9409, 0.5267),  \smallskip\\ 
\tilde\pv_3 = (0.3236, 1.184, 0.3190),
\end{array}\right.
\ea
obtained via least-squares optimization using $\Lambda = [0, 15.89]$, $\Nc_\lambda = 400$, $N_t = 80$ and $T_f = 4.210$. The results are shown in Figure~\ref{fig:multiFreq3SquareNp1ImpEx5}.  The FPI converges, as confirmed by the theory, at a rate faster than the previous case with $\Np=1$ as expected, and we observe that GMRES also converges more rapidly. An examination of the effective convergence
rates, however, shows that when using GMRES there is not much difference in total computational cost between using $\Np=1$ or $\Np=2$.
GMRES converges more rapidly with $\Np=2$ at about twice the CPU cost per iteration, but using about half the memory to store the Krylov vectors.

{
	\newcommand{\figw}{6cm}
	\newcommand{\figh}{5.2cm}
\begin{figure}[htb]
\begin{center}
	\begin{tikzpicture}
		\useasboundingbox (0,.3) rectangle (13,.95*\figh);  

		\begin{scope}[yshift=0*\figh]
			\figByWidth{        0}{0}{fig/sq128originalFD2TSINf2Np1BadCase}{\figw}[0.][0.][0.][0.]
			\figByWidth{1.1*\figw}{0}{fig/sq128originalFD2TSINf2Np1BadCaseFixPointMuFunctionDiscreteBeta}{\figw}[0.][0.][0.][0.]
		\end{scope}  
		
	\end{tikzpicture}
\end{center}
\caption{Multi-frequency algorithm for $\omega_m=\{5.1,10.1\}$ with $N_{p,m}=\{1,1\}$ and traditional filter parameters.  Implicit time-stepping and fourth-order accuracy on a square128 grid. 
   Left: Behavior of the residuals for the FPI and GMRES iterations. Right: MFWH fixed-point convergence function $|\mu_d(\lambdaTilde(\lambda))|$.  The black circle marks the spurious resonance where an eigenvalue (red cross) has $\mu_d=1$.
}
\label{fig:multiFreq2SquareNp1ImpBadCase}
\end{figure}
}

{
	\newcommand{\figw}{6cm}
	\newcommand{\figh}{5.2cm}
	\begin{figure}[htb]
		\begin{center}
			\begin{tikzpicture}
				\useasboundingbox (0,.3) rectangle (13,.95*\figh);  
				
				\begin{scope}[yshift=0*\figh]
					\figByWidth{        0}{0}{fig/sq128leastsquareFD22TSINf2Np1BadCase}{\figw}[0.][0.][0.][0.]
					\figByWidth{1.1*\figw}{0}{fig/sq128leastsquareFD22TSINf2Np1BadCaseFixPointMuFunctionDiscreteBeta}{\figw}[0.][0.][0.][0.]
				\end{scope}  
				
			\end{tikzpicture}
		\end{center}
		\caption{Multi-frequency algorithm for $\omega_m=\{5.1,10.1\}$ with $N_{p,m}=\{1,1\}$ and optimized filter parameters.  Implicit time-stepping and fourth-order accuracy on a square128 grid. 
   Left: Behavior of the residuals for the FPI and GMRES iterations. Right: MFWH fixed-point convergence function $|\mu_d(\lambdaTilde(\lambda))|$.  There are no spurious resonances and both the FPI and GMRES converge to the correct solution for each value of $\omega_m$.
		}
		\label{fig:multiFreq2SquareNp1ImpBadCaseLsq}
	\end{figure}
}

In Figure~\ref{fig:multiFreq2SquareNp1ImpBadCase}, we consider a two-frequency case for which there is a spurious resonance, i.e.~there exists an eigenvalue such that $\mu_d(\lambdaTilde_{h,\nu}) = 1$, $\lambdaTilde_{h,\nu}\neq\omega_m$, $m=1$ or~$2$. This case is constructed by adjusting the domain length $L$ of the square domain $\Omega = [0, L] \times [0, L]$, so that a discrete eigenvalue satisfies $\mu_d = 1$ exactly. It is possible to make this construction because the filter function uses the traditional parameters and, as shown in  the right graph, there are two intervals for which $\vert\mu_d\vert>1$.  We adjust $L$ so that one eigenvalue sits at the right edge of the right-most interval as noted by the black circle in the plot. The Krylov matrix $A_h$ in~\eqref{eq:krylovSystem} becomes singular, and although the Krylov solver converges nicely (left graph) it converges to incorrect solutions for each $m$.  This assessment is made by computing a relative error in the discrete solutions and a relative residual in the discrete equations using
\ba
\label{eq:errorAndResid}
e_h^{(m)} \eqdef \f{\Vert V_\jv^{(m,k_\infty)} - U^{(m)}_{\jv}\Vert}{\|U_\jv^{(m)}\|}, \qquad r_h^{(m)}\eqdef \f{\|f_\jv^{(m)}-\bigl(L_{p,h}-\omega_m^2I\bigr)V_\jv^{(m,k_\infty)}\|}{\|f_\jv^{(m)}\|},
\ea
where $k_\infty$ is an iterate sufficiently close to convergence and the norms are defined here to be
\ba
\Vert Q_\jv\Vert\eqdef \max_{\jv\in\Omega_h}\vert Q_\jv\vert,
\ea
for a generic grid function $Q_\jv$.  Note that the errors and residuals listed in the top row of Table~\ref{tab:gmresIter} are all approximately $10\sp{-2}$, indicating convergence to incorrect solutions.
Figure~\ref{fig:multiFreq2SquareNp1ImpBadCaseLsq} shows results of the MFWH\ algorithm for the same problem, but using optimized parameters given by
\ba
\left\{ \begin{array}{l}
\tilde\pv_1 = (1.046, 0.1588, 1.826), \smallskip\\
\tilde\pv_2  = (-1.224, 1.226, 0.7214),
\end{array}\right.
\ea
for $\Lambda = [0, 28.04]$, $\Nc_\lambda = 120$, $N_t=20$ and $T_f=1.258$.
As shown in the right plot, there are no intervals with $\vert\mu_d\vert>1$ and so no spurious resonances occur.  The Krylov matrix is nonsingular, and both the FPI and GMRES converge to correct solutions for each $m$ as confirmed by the relative errors and residuals reported in the bottom row of~\cref{tab:gmresIter}.

\begin{table}[H]\footnotesize
 	\begin{center}
 		\renewcommand{\arraystretch}{1.5}   
 		\begin{tabular}{|c|c|c|c|c|}
 			\hline
 			Filter & $e_h^{(1)}$ & $e_h^{(2)}$
 			& $r_h^{(1)}$ & $r_h^{(2)}$ \\
 			\hline
 			Traditional&5.36e-02 & 4.25e-02 & 4.18e-02 & 1.35e-02 \\
 			Optimized&9.90e-12 & 7.78e-12 & 8.95e-10 & 1.66e-10 \\
 			\hline
 		\end{tabular}
 	\end{center}
 	\caption{Relative errors $e_h^{(m)}$ and relative residuals $r_h^{(m)}$ for converged solutions of GMRES applied to the MFWH\ algorithm using traditional and optimized filter parameters. The algorithm with traditional filter parameters has a spurious resonance leading to a singular Krylov matrix and convergence to an invalid result.  The algorithm with optimized filter parameters has no spurious resonances.  The Krylov matrix is nonsingular for this case and GMRES converges to the correct discrete solutions of the Helmholtz problems with their respective frequencies.}
 	\label{tab:gmresIter}
 \end{table}
{
	\newcommand{\figw}{6cm}
	\newcommand{\figh}{5.2cm}
	\begin{figure}[htb]
		\begin{center}
			\begin{tikzpicture}
				\useasboundingbox (0,.3) rectangle (13,.95*\figh);  
				
				\begin{scope}[yshift=0*\figh]
					\figByWidth{        0}{0}{fig/sq256leastsquareFD24TSINf7Np4}{\figw}[0.][0.][0.][0.]
					\figByWidth{1.1*\figw}{0}{fig/sq256leastsquareFD24TSINf7Np4FixPointMuFunctionDiscreteBeta}{\figw}[0.][0.][0.][0.]
				\end{scope}  
				
				
			\end{tikzpicture}
		\end{center}
		\caption{Multi-frequency algorithm for $\omega_m = \{ 15, 21, 26, 32, 41, 49, 58 \}$ with $N_{p,m} = [4, 5, 6, 8, 10, 12, 14]$ and optimized filter parameters.  Implicit time-stepping and fourth-order accuracy on a square256 grid. 
   Left: Behavior of the residuals for the FPI and GMRES iterations. Right: MFWH fixed-point convergence function $|\mu_d(\lambdaTilde(\lambda))|$.
		}
		\label{fig:multiFreq7SquareNp4Imp}
	\end{figure}
}

As a final example, Figure~\cref{fig:multiFreq7SquareNp4Imp} shows results for a seven-frequency case defined on a unit square domain
with homogeneous Dirichlet boundary conditions.  The frequencies are taken to be
\bse
\ba
\omega_m = \{ 15, 21, 26, 32, 41, 49, 58 \},
\ea
while the parameters that define the seven Gaussian forcing functions are given by
\ba
\left\{ \begin{array}{l}
   a_{g,m} =   \omega_m^2, \smallskip\\
   b_{g,m} =  \omega_m,  \smallskip \\
  \xv_{0,m} = \{ [.6,.45], [.4,.55], [.55,.5] , [.5,.5] , [.44,.54] , [.53,.45] , [.44,.47]   \} .
\end{array}\right.
\ea
\ese
The optimized filter parameters used here, obtained from the least-squares optimization with $\Lambda = [0, 80.60]$ , $\Nc_\lambda = 840$, $N_t=155$, and $T_f=1.6848$, are $\tilde\pv_m=(\tilde a_m,\tilde c_m,\tilde s_m)$, where
\ba 
\left\{ \begin{array}{l}
	\tilde a_m = \{-4.424, 0.6839, -1.087, 11.90, 0.1628, -32.26, 24.95\}, \smallskip\\
	\tilde c_m = \{1.199, 0.8205, 0.8650, 0.9483, 0.9781, 1.101, 1.078 \}, \smallskip\\
	\tilde s_m = \{0.1287, -0.2387, 0.1575, -0.2345, 0.06153, 0.09895, 0.2539 \} .
\end{array}\right.
\ea 
Although the FPI convergence rate is close to one, GMRES converges nicely.  A base value $\Np=4$ is used for these calculations, which in turn implies $N_{p,m} = [4, 5, 6, 8, 10, 12, 14]$ for all frequencies.  From the right graph of Figure~\cref{fig:multiFreq7SquareNp4Imp} we observe that since larger frequencies have smaller periods, their corresponding value for $N_{p,m}$ is larger and thus the widths of the main peaks near each $\omega_m$ are roughly the same for large and small frequencies.

Figure~\cref{fig:multiFreq7SquareNp4ImpSol} shows color contours of the computed MFWH solutions $V_{\jv}^{(5,k_\infty)}$ and $V_{\jv}^{(7,k_\infty)}$ corresponding to the frequencies $\omega_5 = 41$ and $\omega_7 = 58$, respectively, taken at the final GMRES iterate $k_\infty=165$ without restart.  Each solution is forced by a Gaussian forcing centered at $\xv_{0,5} = [0.44,0.54]$ and $\xv_{0,7} = [0.44,0.47]$, respectively, and indicated by the black dots in the plots. Also shown is the combined solution, i.e.~$\sum_{m=1}^{\Nf} V_{\jv}^{(m,k)}$. The individual WaveHoltz solutions $V_{\jv}^{(m,k)}$, $m=1,2,\ldots,\Nf$, are recovered from this 
combined solution by the MFWH time-filtering procedure.

{
	\newcommand{\figw}{4.5cm}
	\newcommand{\figh}{4.4cm}
	\begin{figure}[htb]
		\begin{center}
			\begin{tikzpicture}
				\useasboundingbox (0.25,.35) rectangle (3*\figh,4.6); 

				\begin{scope}[yshift=0*\figh]
					\figByWidth{        0*\figh}{0}{fig/sq256leastsquareFD4TSINf7Np4SolutionSumPlot}{\figw}[0.][0.][0.][0.]
					\figByWidth{1*\figh}{0}{fig/sq256leastsquareFD4TSINf7Np4Solution5Plot}{\figw}[0.][0.][0.][0.]
					\figByWidth{2*\figh}{0}{fig/sq256leastsquareFD4TSINf7Np4Solution7Plot}{\figw}[0.][0.][0.][0.]
				\end{scope}
				\draw[->, thick, black] (1.80, 0.45) -- (2.30, 2.55);
   			\node[above] at (1.80, 0.06) {\tiny forcing};
   			\draw[->, thick, black] (6.00, 0.45) -- (6.52, 2.70);
   			\node[above] at (6.00, 0.06) {\tiny forcing};
   			\draw[->, thick, black] (10.40, 0.45) -- (10.95, 2.50);
   			\node[above] at (10.40, 0.06) {\tiny forcing};
			\end{tikzpicture}
		\end{center}
		\caption{Color contours of solutions of the MFWH algorithm for a seven-frequency case: $\omega_m = \{ 15, 21, 26, 32, 41, 49, 58 \}$ with $N_{p,m} = [4, 5, 6, 8, 10, 12, 14]$ and optimized filter parameters.  Implicit time-stepping and fourth-order accuracy on a square256 grid. Left: Sum of all seven solutions at convergence. Middle: Solution for $\omega_5=41$. Right: Solution for $\omega_7=58$.  The black dots in the three plots show the locations of the relevant Gaussian sources.}
		\label{fig:multiFreq7SquareNp4ImpSol}
	\end{figure}
}

  In~\cite{overHoltzPartOne} the $O(N)$ scaling of the WaveHoltz algorithm for a single frequency was shown. 
A key ingredient was the use of implicit time-stepping with a fixed number of steps per period, together with an $O(N)$ solver for the implicit time-stepping equations.
Figure~\ref{fig:squareGridCpuScaling} shows MFWH scaling results when solving for three solutions using the same parameters as in~\cref{eq:Nf3parameters}.
The implicit time-stepping equations are solved with the Ogmg multigrid solver~\cite{OGMG,automg,multigridWithNonstandardCoarsening2023} to a tolerance of $10^{-10}$.
The left table in Figure~\ref{fig:squareGridCpuScaling} shows that the number of wave-solves and average number of multigrid cycles per wave-solve remain
fairly constant as the grid is refined. A ``wave-solve'' refers to a single MFWH iteration consisting of solving the wave equation followed by the time filtering.
The column titled ``CPU ratio'' shows the ratio of CPU time of the current row to that of the previous row. 
This ratio should be about $4$ since the number of grid points $N$ increases by a factor of $4$ from one row to the next.
The right graph plots the relative CPU times scaled by the number of grid points. After an initial decrease the curves are relatively flat consistent with an $O(N)$ scaling.

\newcommand{\mfwhMGTable}{%
\footnotesize
\begin{tabular}{|c|c|c|c|c|c|} \hline
      \multicolumn{6}{|c|}{MFWH scaling: square, order 4} \\ \hline 
               &   grid-       &  wave-    &  multigrid    &  CPU          &  CPU      \\
 $\dx $       &   points      &  solves   &  cycles       &  seconds      &  ratio    \\ \hline
 $1/64   $    &   4.2e+03     &  $22$   &  $7.1$    &   $  5.6$      &            \\
 $1/128  $    &   1.7e+04     &  $22$   &  $7.0$    &   $ 14.3$      &      2.55       \\
 $1/256  $    &   6.6e+04     &  $22$   &  $7.0$    &   $ 49.9$      &      3.50       \\
 $1/512  $    &   2.6e+05     &  $20$   &  $6.7$    &   $169.8$      &      3.40       \\
 $1/1024 $    &   1.1e+06     &  $20$   &  $6.7$    &   $690.9$      &      4.07       \\
 $1/2048 $    &   4.2e+06     &  $20$   &  $6.7$    &   $2862.0$      &      4.14       \\
\hline
\end{tabular}    
}%

{
\newcommand{\figWidth}{5.5cm}
\newcommand{\figHeight}{4.5cm}

\begin{figure}[htb]
\begin{center}
\resizebox{\textwidth}{!}{
\begin{tikzpicture}[scale=1]
  \useasboundingbox (0,.4) rectangle (2.6*\figWidth,1.05*\figHeight);  

   \begin{scope}[xshift=4cm,yshift=2.5cm]
     \node  { \mfwhMGTable };
    \end{scope}
   \begin{scope}[xshift=1.6*\figWidth]
    \figByWidth{0}{0}{fig/squareMFWHCpu}{\figWidth}[0][0][0][0];
  \end{scope}  

\end{tikzpicture}
}
\end{center}
\caption{Solving a three frequency MFWH problem with implicit time-stepping and multigrid.
  Left: Solver statistics for grids of increasing resolution.
  Right:  Normalized values of CPU$(N)$/$N$, versus number of grid points.
   The results show the near optimal $O(N)$ scaling of the MFWH algorithm.
    }
\label{fig:squareGridCpuScaling}
\end{figure}
}

\section{Conclusions} \label{sec:conclusions}

A new multi-frequency WaveHoltz (MFWH) 
algorithm has been presented to simultaneously compute multiple Helmholtz solutions for different frequencies and different forcing functions.
The MFWH algorithm is an extension of the original, single-frequency WaveHoltz algorithm.  The new algorithm is based on the solution of a single wave equation that is forced
by a composite forcing function that is the sum of different component functions that oscillate at different frequencies.
At each stage in the iteration multiple time filters are applied to the solution of the wave equation in order to update the approximations
to the different Helmholtz solutions.
The MFWH algorithm was shown to be a fixed-point iteration whose convergence behavior depended on the frequencies of the Helmholtz problems and on the parameters of the filter function.
The convergence behavior was analyzed using eigenfunction expansions and the asymptotic convergence rate was found.
Using the filter parameters corresponding to the original WaveHoltz algorithm, it was found that the MFWH fixed-point iteration may not converge due to the presence of intervals for which $\vert\mu\vert>1$.  A least-squares optimization approach was described to select filter parameters so that $\vert\mu\vert\le1$ and the MFWH fixed-point iteration converged (for all problems tested).
Moreover, it was shown that the iteration could be accelerated using a Krylov-based scheme such as GMRES.  With the least-squares parameters, spurious resonances were avoided and GMRES converged.
Numerical results in two dimensions confirmed the theory and demonstrated the $O(N)$ scaling of the algorithm.

\appendix

\section{Matlab codes to determine optimal parameters} \label{sec:matlab_codes_to_determine_optimal_parameters}

Matlab codes that determine the optimal filter parameters, as described in Section~\ref{sec:optimization}, are given in Listing~\ref{code:optparams}. Algorithm~\ref{alg:SSWHV2} is formulated in terms of the continuous filter functions in~\eqref{eq:compFuncs}, while the Matlab codes employs numerical quadratures based on the grid used for the time integration of the wave equation in the MFWH algorithm. Thus, in addition to the vector of frequencies $[\omega_1,\ldots,\omega_\Nf]$ and $N_{p,1}$, the input to the code takes the total number of equally-spaced time steps $N_t$ in the integration to the final $T_{f,1}=2\pi N_{p,1}/\omega_1$.  A standard trapezoidal rule approximation is used for the interval $[0,T_{f,1}]$, while the quadrature weights given in~\eqref{eq:DefNI} and~\eqref{eq:betadDefA} are used for the integrations over the subintervals $[0,T_{f,m}]$, $m=2,\ldots,\Nf$, which are not aligned with the base grid in general.  Note also, that the modified frequencies $\tilde\omega_m$ in~\eqref{eq:omegaImplicit} with corresponding modified periods are used as input to the code for the results presented in Section~\ref{sec:numerical_confirmation_of_the_convergence_analysis}.

\begin{lstlisting}[language=matlab, numbers=left, stepnumber=1, firstline=1,caption={Code for optimized parameters.},label=code:optparams,frame=single]
function [xa_opt, xs_opt, xc_opt] = OptParametersC(freqArray,Np,Nt,Ndelta,cdelta)
% OptParametersAll: Computes optimal MFWH filter parameters via least-squares.
%
% Solves the linearly constrained least-squares problem to find the optimal
% filter parameters that minimize the multi-frequency WaveHoltz filter
% function mu(lambda) away from the target frequencies omega_m, subject to
% the interpolation and flatness constraints at each omega_m.
%
% The target frequencies are ordered increasingly: omega_1 < ... < omega_Nf
% INPUT:
%     freqArray : target frequencies [omega_1,...,omega_Nf]
%     Np        : number of periods for T1
%     Nt        : number of time steps
%     Ndelta    : points per half-width
%     cdelta    : safety factor
%     
%
% OUTPUT:
%   xa_opt : optimal a-parameters (Nf x 1)
%   xs_opt : optimal s-parameters (Nf x 1)
%   xc_opt : optimal c-parameters (Nf x 1)
%
% METHOD:
%   1. Assemble constraint matrix G = [Fc Fs; Fc' Fs'] and H = [Fa; Fa']
%      evaluated at the target frequencies omega_m.
%   2. Select Nlambda midpoint quadrature points over [0, lamMax].
%   3. Evaluate basis functions fc, fs, fa at quadrature points.
%   4. Form and solve the unconstrained least-squares problem for xa_opt.
%   5. Recover xc_opt and xs_opt from the constraint system.
%
% See also: TrapRule, TrapWeights

%% Time parameters
freqArray = freqArray(:);          %column vector

T_m  = 2*pi ./ freqArray;
N_pm = floor(Np * (T_m(1)./T_m));  % number of periods per frequency
T_fm = N_pm .* T_m;                % final times per frequency
nf   = length(freqArray);          % number of frequencies
Tf   = max(T_fm);
dt   = Tf/Nt;
                
%% Precompute quadrature weights (done once, reused in all filter evaluations)
sigma = TrapWeights(Nt, nf, T_fm, dt);

%% Time integration operator (trapezoid rule with frequency-specific weights)
% Int(g, lambda, om, freq) = (2/T_freq) * integral_0^{T_freq} sigma*g dt
Int = @(g,lambda,om,freq) 2/T_fm(freq)*TrapRule(g,lambda,om,freq,Nt,dt,sigma);

%% Basis filter functions evaluated at a vector x, for frequency index m
% f_c(x;omega_m) = (2/T_fm) int_0^{T_fm} cos(omega_m t) cos(lambda t) dt
% f_s(x;omega_m) = (2/T_fm) int_0^{T_fm} sin(omega_m t) cos(lambda t) dt
% f_a(x;omega_m) = (2/T_fm) int_0^{T_fm}                cos(lambda t) dt
% and their derivatives with respect to lambda (denoted with prime):
fc  = @(x,m) Int(@(omega,lambda,t)  cos(omega*t).*cos(lambda*t), x,freqArray(m),m);
fs  = @(x,m) Int(@(omega,lambda,t)  sin(omega*t).*cos(lambda*t), x,freqArray(m),m);
fa  = @(x,m) Int(@(omega,lambda,t)  cos(lambda*t),               x,freqArray(m),m);
fcp = @(x,m) Int(@(omega,lambda,t) -t.*cos(omega*t).*sin(lambda*t), x,freqArray(m),m);
fsp = @(x,m) Int(@(omega,lambda,t) -t.*sin(omega*t).*sin(lambda*t), x,freqArray(m),m);
fap = @(x,m) Int(@(omega,lambda,t) -t.*sin(lambda*t),               x,freqArray(m),m);

%% Assemble constraint matrices G, H and right-hand side z
% G = [Fc  Fs ] in R^{2Nf x 2Nf},  H = [Fa ] in R^{2Nf x Nf},  z = [1;0]
%     [Fc' Fs']                        [Fa']
% where G(m,n) = fc(omega_m; omega_n), G(m,n+nf) = fs(omega_m; omega_n), etc.

G = zeros(2*nf, 2*nf);
H = zeros(2*nf, nf);
z = zeros(2*nf, 1);
for m = 1:nf
    for n = 1:nf
        G(m,n)       = fc(freqArray(m), n);    % Fc block
        G(m,n+nf)    = fs(freqArray(m), n);    % Fs block
        G(m+nf,n)    = fcp(freqArray(m), n);   % Fc' block
        G(m+nf,n+nf) = fsp(freqArray(m), n);   % Fs' block
        H(m,n)       = fa(freqArray(m), n);    % Fa block
        H(m+nf,n)    = fap(freqArray(m), n);   % Fa' block
    end
    z(m) = 1;   % interpolation constraint: mu(omega_m) = 1
end

%% Midpoint quadrature points over Lambda = [0, lamMax]
% lamMax extends slightly beyond omega_Nf by cdelta*nf half-widths delta_Nf
% Nlambda chosen so that each half-width contains Ndelta quadrature points
%Ndelta   = 20;                                   % points per half-width
%cdelta   = 2;                                    % safety factor
omegaNf  = max(freqArray);                        % largest frequency
deltaNf  = omegaNf / (2*N_pm(end));               % half-width of last peak
lamMax   = omegaNf + cdelta*nf*deltaNf;           % upper limit of Lambda
Nlambda  = ceil(Ndelta * lamMax / deltaNf);       % number of quadrature points
lambdaVals = ((1:Nlambda) - 1/2) * lamMax/Nlambda;  % midpoint rule points

%% Evaluate basis functions at quadrature points (each nf x Nlambda)
% Vectorfc(:,k) = [fc(lambda_k;omega_1); ...; fc(lambda_k;omega_Nf)]
Vectorfc = zeros(nf, Nlambda);
Vectorfs = zeros(nf, Nlambda);
Vectorfa = zeros(nf, Nlambda);
for k = 1:Nlambda
    for n = 1:nf
        Vectorfc(n,k) = fc(lambdaVals(k), n);
        Vectorfs(n,k) = fs(lambdaVals(k), n);
        Vectorfa(n,k) = fa(lambdaVals(k), n);
    end
end

%% Form and solve the unconstrained least-squares problem for xa
% After eliminating xcs = G\(z - H*xa) from the constraint G*xcs + H*xa = z,
% mu(lambda) = fcs^T * G\z + (fa^T - fcs^T * G\H) * xa = mu0 + ga^T * xa
% Minimize int_Lambda |mu(lambda)|^2 dlambda => M*xa = m0
Vectorfcs = [Vectorfc; Vectorfs];           % (2nf x Nlambda)
% GinvZH    = G \ [z H];                    % factorize G once
% GinvZ     = GinvZH(:,1);                  % G\z
% GinvH     = GinvZH(:,2:end);              % G\H
M         = Vectorfa'  - Vectorfcs'*(G\H);  % least-squares matrix (Nlambda x nf)
m0        = -Vectorfcs'*(G\z);              % least-squares rhs   (Nlambda x 1)

fprintf('Optimized filter via least squares\n');
xa_opt = M \ m0;                            % solve for xa (nf x 1)

%% Recover xc and xs from the linear constraint
y_opt  = G\(z - H*xa_opt);             % [xc; xs] = G\(z - H*xa)
xc_opt = y_opt(1:nf);                      % c-parameters (nf x 1)
xs_opt = y_opt(nf+1:end);                  % s-parameters (nf x 1)

end

% =========================================================================
function f = TrapRule(g, lambda, om, freq, Nt, dt, sigma)
% TrapRule: Evaluates the weighted time integral of a filter kernel g.
%
% Computes f(i) = sum_{n=0}^{N} sigma(n+1,freq) * g(om, lambda(i), t_n)
% for each entry lambda(i) of the input vector lambda, where t_n = n*dt.
%
% INPUT:
%   g      : kernel function handle g(omega, lambda, t)
%   lambda : vector of evaluation points (1 x nx)
%   om     : frequency value omega_m (scalar)
%   freq   : frequency index m (used to select sigma column)
%   Nt     : number of time steps
%   dt     : time step size
%   sigma  : quadrature weights (nStep+1 x nf), from TrapWeights
%
% OUTPUT:
%   f : vector of integral values (1 x nx)

nx      = length(lambda);
f       = zeros(1, nx);
for i = 1:nx
    s = 0;
    for n = 1:Nt+1
        t = (n-1)*dt;
        s = s + sigma(n,freq)*g(om, lambda(i), t);
    end
    f(i) = s;
end
end

% =========================================================================
function sigma = TrapWeights(Nt,nf,T_fm, dt)
% TrapWeights: Computes frequency-specific quadrature weights for TrapRule.
%
% Uses a modified trapezoid rule when the frequency-specific final time
% T_fm(j) does not coincide with a point on the common time grid.
% the integration interval [0, T_fm] may not align exactly with a time step,
% so a smoothly varying weight is applied near the endpoint T_fm.
%
% INPUT:  
%     dt   : time step size
%     Nt   : number of time steps
%     nf   : number of frequencies Nf
%     T_fm : final integration times [T_f1,...,T_fNf]
%
% OUTPUT:
%   sigma : quadrature weight matrix (nStep+1 x nf)
%           sigma(n,j) = weight for time step n and frequency j

sigma = zeros(Nt+1, nf);
k     = zeros(nf, 1);

for j = 1:nf
    k(j)  = floor(T_fm(j)/dt);                    % last full time step index
    beta  = (dt*(k(j)+1) - T_fm(j))/dt;           % fractional overshoot

    sigma(1,j)     = 1/2*dt;                      % left endpoint weight
    for n = 1:k(j)-1
        sigma(n+1,j) = dt;                        % interior weights
    end
    sigma(k(j)+1,j) = 1/2*dt*(2 - beta^2);       % near-endpoint weight
    if k(j) < Nt
        sigma(k(j)+2,j) = 1/2*dt*(1 - beta)^2;   % final partial step weight
    end
end
end
\end{lstlisting}

\bibliographystyle{siam}
\bibliography{journal-ISI,henshaw,henshawPapers,waveHoltz,helmholtz,appelo}

\end{document}